\documentclass{article}
\usepackage{latexsym}
\usepackage{graphicx}
\usepackage{amssymb}

\def\A{{\mathcal A}}
\def\Aut{\hbox{Aut}}

\def\L{{\mathcal L}}

\def\RR{{\mathbb{R}}}
\def\Id{{\mathbb{I}}}
\def\W{{\mathcal W}}
\def\Prob{\hbox{Prob}}

\def\E{\hbox{E}}
\def\Var{\hbox{Var}}
\def\Cov{\hbox{Cov}}

\def\bfu{{\underline{u}}}

\newtheorem{theorem}{Theorem}[section]
\newtheorem{lemma}[theorem]{Lemma}

\newtheorem{corollary}[theorem]{Corollary}

\newenvironment{proof}[1][Proof]{\begin{trivlist}
\item[\hskip \labelsep {\bfseries #1}]}{\end{trivlist}}
\newenvironment{definition}[1][Definition]{\begin{trivlist}
\item[\hskip \labelsep {\bfseries #1}]}{\end{trivlist}}

\begin{document}

\title{Multiple pattern matching: A Markov chain approach}
\author{Manuel E. Lladser  \and M. D. Betterton \and Rob Knight
\thanks{We are grateful to Robert S. Maier and several researchers in the Analysis of Algorithms (AofA) community for their helpful comments and suggestions while preparing this manuscript.}
}

\maketitle

\begin{abstract}
  RNA motifs typically consist of short, modular patterns that include
  base pairs formed within and between modules. Estimating the
  abundance of these patterns is of fundamental importance for
  assessing the statistical significance of matches in genomewide
  searches, and for predicting whether a given function has evolved
  many times in different species or arose from a single common
  ancestor. In this manuscript, we review in an integrated and
  self-contained manner some basic concepts of automata theory,
  generating functions and transfer matrix methods that are relevant
  to pattern analysis in biological sequences. We formalize, in a
  general framework, the concept of Markov chain embedding to analyze
  patterns in random strings produced by a memoryless source. This
  conceptualization, together with the capability of automata to
  recognize complicated patterns, allows a systematic analysis of
  problems related to the occurrence and frequency of patterns in
  random strings. The applications we present focus on the concept of
  synchronization of automata, as well as automata used to search for
  a finite number of keywords (including sets of patterns generated
  according to base pairing rules) in a general text.
\end{abstract}

\section{Introduction}
\label{intro}

The importance of RNA in biology is increasing as we learn more about
the function of RNA molecules. Some RNA molecules are passive
messengers in translation (a step in the production of protein
molecules from the DNA genome), but RNA molecules can also act as a
catalysts \cite{cech81,guer83}. Recent estimates suggest that the
human genome may encode up to 75,000 small RNA genes, which is at
least three times the number of protein-coding genes \cite{lu05}.
Because new functional RNA molecules are being discovered every day,
the problem of understanding the structure and sequence requirements
for RNA function is of increasing importance.

Functional RNA molecules share important structural and sequence
characteristics. These RNA molecules typically consist of short,
evolutionarily conserved regions (modules) that are separated by
essentially random spacer sequences that can vary both in length and
nucleotide sequence \cite{knig03}. Modules often base pair with each
other, an effect which introduces long-range correlations among parts
of the sequence. 
(For more detailed definitions of patterns and modules, see below.)

If a particular motif corresponds to a functional RNA molecule, the
corresponding modular pattern may be statistically over- or
underrepresented in the genome. This assumption is used in genomewide
searches for possible functional RNA molecules. The estimation of
over- or underrepresentation
requires us to calculate the probability that the modular pattern
occurs in some statistical model of the genome sequence. Therefore, the
study of RNA sequences is directly related to pattern matching and the
probability of occurrence of patterns in random strings.

Traditionally, sequence similarity between genetic sequences in
different organisms has been interpreted to mean that the gene in both
organisms share a common ancestor. This assumption underlies many
sequence analysis algorithms. However, increasing evidence suggests
that sequence similarity may not always imply common descent of RNA molecules.
This may occur because, despite the diversity of functional RNA
molecules, some functions can only be evolved in a relatively small
number of ways. For example, the hammerhead ribozyme, a self-cleaving
RNA that has an evolutionarily conserved catalytic core of only 11
nucleotides, has both been observed in a wide range of organisms and
has also been artificially selected from random-sequence backgrounds
\cite{tang00,sale01}. Similarly, artificial selection of RNAs from
random sequences has recaptured the sequence of the catalytic core of
the ribosome \cite{welc97,yaru00} and features of the genetic code
\cite{yaru05}.  Therefore, study of RNA molecules may require new
models of sequence evolution that characterize the origins of a motif
from a random sequence.

Although probabilistic models for sequence evolution from a common
ancestor are well-established \cite{kimu81,fels81}, probabilistic
models for independent origins of an RNA motif in random-sequence
backgrounds have been less well studied \cite{knig03,knig05}.  The
long-range correlations introduced by base pairing can be difficult to
accommodate in search algorithms.  Paired RNA motifs cannot be
represented as regular languages, but instead must be represented as
context-free grammars for full generality \cite{eddy94,riva00}.
Genomewide searches have been performed for several functional RNA
motifs \cite{bour99,ferb00,klei03,grif05}. In these searches, the
statistical significance of matches has typically been assessed using
Monte Carlo simulations, in which the search is repeated using
randomized versions of the search text. This procedure has significant
limitations: it is time-consuming and cannot accurately estimate the
low $p$-values that are important for computing likelihood ratios for
rare events. Therefore, new methods for computing the probability of
sequences in random strings is important for determining the
statistical significance of RNA motif searches.

Previous work has focused on patterns related to RNA structure.
However, recent work has developed other pattern-matching problems
related to RNA sequences.  For example, multiple short protein- or
RNA-binding sequence motifs can combine to regulate a range of
biological processes, including splicing and polyadenylation
\cite{sing06}.  Similarly, 6-base seed sequences that bind short microRNA
molecules (miRNAs) appear to work in concert to repress
translation \cite{lewi05}. Current evidence suggests that while the
motifs function combinatorially (several motifs must be present
together for biological function), no results suggest that specific
base pairing between the modules is required for function. Therefore,
for sequence analysis the motifs can be treated as uncorrelated
(although the sequences within each module can be compound).  

RNA and other biological sequences are not intrinsically random.
However, computational biologists model these sequences as random
(using different models) to assess which patterns within the
sequence are likely to be biologically significant. Therefore the
modeling of RNA as a random sequence is a mathematical construction
rather than a biophysical model of RNA.

Here we review a broad selection of approaches to estimating the
expected number of matches to, or the probability of occurrence of,
RNA motifs. We consider motifs both with and without correlations
(such as base pairing). These approaches draw from many branches of
mathematics and computer science.  Progress in this field has been
limited by the difficulty of integrating results from different fields
that use different concepts and terminology (see section
\ref{sec:glossary} for a glossary of terms). Terminology and previous
work is reviewed in section \ref{sec:previous}. In section
\ref{sec:automata}, we give mathematical definitions and proofs of key
concepts in deterministic pattern matching. These concepts include the
use of automata to search for keywords in databases, and
synchronization of automata to search for multiple patterns
simultaneously. We also give independent proofs that the Aho-Corasick
automaton matches compound patterns, even in cases in which keywords
are subpatterns of other keywords.  In section \ref{sec:markov} we
formalize the Markov chain embedding technique for probablistic
pattern matching.

Sections \ref{sec:app1} and \ref{sec:app2} describe examples where the
methods are applied to sequence analysis problems relevant to RNA
motif searches. The examples we give here, which rely on memoryless
sources, demonstrate how automata theory is a fundamental tool to
analyze the occurrence of patterns in random strings. We also provide references that
extend these examples to Markovian models, and to other, more complex,
models. Our examples rely on the matrix representation of
probabilities and generating functions extracted from the graphs
associated with the automata. As expected, the generating functions are rational functions (i.e., ratios of polynomials). This
feature allows the use of well-known techniques to analyze their
asymptotic behavior for long random strings, such as those encountered
in genomewide analysis.  We will focus primarily on sooner-times
(i.e., the first occurrence of any item from a list of patterns in a
random sequence) and count statistics (i.e., the number of occurrences of a pattern in a random sequence).

\section{Glossary of terms} 
\label{sec:glossary}

\begin{description}
\item[alphabet:] a finite set of characters used to build a text. Example: for
  RNA sequences, the alphabet contains the four nucleotides $A$, $C$,
  $G$, and $U$.

\item [autocorrelation polynomial:] a polynomial in one variable that
  quantifies the degree to which a word overlaps with itself.
  
\item[automaton:] see deterministic finite automaton; also called state machine.
  
\item[Bernoulli source:] a model for the generation of a random string
  in which the probability of a given character is fixed, independently
  of the characters appearing elsewhere in the string; also called a memoryless
  model.

\item[character:] an element of an alphabet; also called letter or symbol. 
 
\item[compound pattern:] any finite set of strings, usually but not
  always consisting of words with a common or similar structure; also
  called a degenerate pattern. Example: $AAC\{U,T\}CCG$ is a compound pattern of two 7-letter strings where the fourth letter in either string can be either $U$ or $T$.
  
\item[correlation:] a dependency between two positions in a pattern,
  such as that introduced by modeling base pairing in RNA or DNA. 
  Example: a correlation would exist if positions 3 and 7 in a pattern
  must base pair with each other; these positions can be filled by any
  letters as long as they form a base pair.
  
\item[De Bruijn graph:] an automaton, the states of which track the
  last $k$ characters read in a text.
  
\item[deterministic finite automaton:] an abstract representation of a (regular pattern) search algorithm, with a finite set of states and rules that specify transitions between states. It is usually represented as a graph composed of a finite number of states (the nodes), transitions between states (the edges) and actions performed upon entering or leaving a state (conveyed by labels on the edges); also called finite-state automaton or finite-state machine. For keywords matching, there is an initial state representing the state before any characters have been matched, and final states representing matches with the keywords.
  
\item[dynamic source:] a particular type of probabilistic model for the generation of a
  random sequence, for which the probability of a character may depend
  on all the preceding characters. Bernoulli and Markovian sources are
  particular instances of dynamic sources.

\item[edit distance:] the distance between two strings, as calculated
  by summing the cost of each elementary operation (e.g., character
  insertion, deletion or substitution) required to convert one string
  into the other.
  
\item[generalized word:] a compound pattern for which all words in the
  pattern have the same length.
  
\item[generating function:] a function in one or more variables for which
  the Taylor series coefficents of the function correspond to
  probabilities or expected values associated with a discrete random
  variable or vector of interest.
  
\item[hidden pattern:] a pattern which may appear separated into blocks rather than as a single and continuous block within the text. Example: the pattern $AGA$ appears 
four times in the string $ACAGCCUGA$ as a hidden pattern.

\item[keyword:] see string.

\item[language:] see pattern.
  
\item[letter:] see character.
  
\item[Markov chain:] a sequence of random variables (here taking
  values in an alphabet) for which the probability of the value taken by a
  random variable is determined by the values of the $k$ previous
  variables. The parameter $k\ge0$ is a finite constant that corresponds to the Markov order of the chain.

\item[Markovian source:] a model for the generation of a random
  sequence for which the probability of a character depends on the $k$
  preceding characters, where $k\ge0$ is the Markov order of the
  sequence.
  
\item[memoryless source:] see Bernoulli source.

\item[module:] see modular pattern.

\item[modular pattern:] an ordered list of simple or compound patterns
  that may include correlations within or between patterns. Example:
  in the context of the DNA alphabet, $AA12...2'1'GT$ is a pattern with two correlated modules, namely $AA12$ and $2'1'GT$, where 1, 2 and $1'$, $2'$ denote correlations between the first and second module. Here $A'=T$, $T'=A$, $G'=C$ and $C'=G$.

\item[non-overlap counting:] the total number of substrings of a text
  that match a pattern, as the text is read from left to right, where a
  given substring can only be considered once for a match; also called
  renewal counting.

\item[overlap counting:] the total number of substrings of a text that
  match a pattern.
  
\item[pattern:] a set of strings. The strings in the pattern
  usually (but not necessarily) are similar to each other. These
  include simple patterns, modular patterns, correlated modular
  patterns, and  any set of words specified by a regular
  expression; also called language.
  
\item[prefix:] a substring which corresponds to the start of a string. Every
  word is a prefix of itself. Example: $A$, $AA$, $AAC$, and $AACU$
  are prefixes of the string $AACUCCG$.
  
\item[reduced pattern:] a pattern where no string is a substring of another string in the pattern.

\item[renewal counting:] see non-overlap counting.

\item[regular expression:] a string that describes all and only those strings belonging to a regular language.

\item[regular language:] see regular pattern.

\item[regular pattern:] a set of strings that can be recognized by a deterministic finite automaton; also called regular language.
  
\item[run:] a maximal sequence of identical characters in a text. Example: the binary string 0010011110 has three runs of zeros, namely 00, 00 and 0, and two runs of ones, 1 and 1111, respectively.

\item[state machine:] see deterministic finite automaton; also called automaton.

\item[string:] a specific sequence of alphabet characters; also called keyword or word.

\item[substring:] a consecutive list of characters within a string; also called sub-word.

\item[suffix:] a substring which corresponds to the ending of a string. Every word is a suffix of itself. Example: $G$, $CG$, $CCG$, and $UCCG$ are suffixes of the string $AACUCCG$.

\item[sub-word:] see substring.

\item[symbol:] see character.
  
\item[sooner-time:] the time required before an event is observed;
  here corresponds to the length of text that precedes the first
  occurrence of a pattern.
  
\item[simple pattern:] a pattern where each position in the string is
  exactly specified by one letter. Example: $AACUCCG$ is a 7-letter
  simple pattern. 
  
\item[suffix:] a substring which corresponds to the end of another string.
  Every word is a suffix of itself. Example: $G$, $CG$, $CCG$, and
  $UCCG$ are suffixes of the string $AACUCCG$.
  
\item[text:] a usually long sequence of characters in which patterns may occur. May be
  randomly generated according to a probabilistic model.
  
\item[transition matrix:] the matrix that summarizes the probability
  that a Markov chain undergoes a transition from one state to
  another.

\item[transfer matrix:] a matrix with polynomial entries, here used to keep track of the number of visits that a Markov chain makes to a certain set of states.

\item[word:] see string.

\end{description}

\section{Prior work on pattern matching}
\label{sec:previous}

In this section we summarize much of the work that uses automata and Markov chains to study patterns in random strings. For an introduction to automata theory and regular expressions see~\cite{HopUll79,Sip96}. See~\cite{Wat95,RobRodSch05} for an introduction to pattern analysis of biological sequences. A comprehensive discussion of patterns in random strings can be found in the book of Lothaire et al.~\cite{Lot05}. Other references give useful background on the mathematical techniques discussed in this paper. An introductory treatment of generating function methods can be found in~\cite{Wil94}. See \cite{FlaSed06} for a broader discussion on generating function and transfer matrix methods. Supplementary references on Markov chains include~\cite{Bre98,Dur99}; for a more detailed discussion of Markov chains
the reader is referred to~\cite{Fel68,Dur04}.

\subsection{Terminology}

\subsubsection{Determistic versus probabilistic pattern matching.}

Early work in computer science focused on \textbf{deterministic}
pattern matching, where the text to be analyzed is given (for example,
the abstract of this paper) and one wishes to search the text for a given pattern.  The
number of occurrences of the pattern has a definite answer once the
text is given.  For applications of pattern matching to problems in
biological sequence data, one is typically interested in probabilistic
versions of the pattern matching problem. Therefore this review
focuses on research in \textbf{probabilistic} pattern matching. Here
one models the biological sequence as a random string produced
according to some model. Typically, the sequence is assumed to be
produced by a memoryless source or a Markovian source, although hidden
Markov models are sometimes used. The pattern matching problem
can then be formulated as a probabilistic question; different papers
answer slightly different questions. Previous work can be categorized
as problems involving (i) counting (what is the probability that
a given pattern occurs $m$ times in a random string of length $n$?),
(ii) occurrence (what is the probability that a given pattern
occurs or does not occur in a random string of length $n$?), (iii) type of occurrence (what is the probability that a string in a pattern is the first one observed?), and (iv) distance between occurrences (what is the typical
distance between successive occurrences of a pattern in a random
string of length $n$?). Note that the question of occurrence
probability (ii) is a special case of the counting problem (i).
  
\subsubsection{Type of pattern.}

Research in this field has considered a range of different kinds of
pattern. The most basic case is a \textbf{simple pattern}. A simple
pattern is a string where each position in the string is exactly
specified by one letter; the word {\it dog} would be an example of a
simple pattern. A \textbf{compound pattern} is a finite set of simple patterns; for example, a keyword search for the words
{\it dog}, {\it cat}, and {\it snake} would seek to match 
a compound pattern. Compound patterns
are sometimes specified by letting some positions in the string 
be chosen from a range of characters. For example, the words {\it snake} and {\it snare} could be represented by the compound pattern {\it sna\{k,r\}e}.  A pattern is referred to as \textbf{forbidden} if the pattern matching problem seeks to exclude
occurrences of the pattern, rather than find occurrences of it.

A \textbf{correlated pattern} contains positions where characters must
be related by some rule. For example, one could search for the
correlated pattern {\it 1o1} with a rule that positions marked by the
number 1 must be the same letter. A search for such a correlated
pattern would find all 3-letter words where the first and third
letters are the same, such as {\it mom} and {\it tot}.

A \textbf{modular pattern} is composed of subpatterns which must appear in a certain order, but which could be separated by one or more characters. For example, a
search for the modules {\it cat...dog} would match examples in the
text where the word {\it cat} occurs, followed by any number of
characters, followed by the word {\it dog}.  The number of characters
allowed between the modules can be unbounded, bounded, or
specified uniquely. For instance, {\it cat\#\#...dog} would match
examples where {\it cat} occurs, followed by at least two characters
followed by {\it dog}. A modular pattern could also contain
correlations within or between modules. A modular pattern 
may include an infinite number of simple patterns.

A \textbf{regular pattern} is a pattern that can be described by a regular expression of the type used in computer science.

\subsubsection{Overlaps.}

When matching more complicated patterns, one must specify how to deal
with overlaps of words. The \textbf{overlapping count} of a pattern in
a text corresponds to the number of substrings of the text that belong
to the set of words specified by the pattern. For example, there are 4
overlapping occurrences of {\it TATA} and 5 overlapping occurrences of
{\it ATA} in the text {\it ATATATATATA}; therefore, there are 9
overlapping occurrences of the compound pattern {\it TATA, ATA} in 
this text. There are only 3 overlapping occurrences of the modular pattern {\it TA\#...TATA} in the text {\it ATATATATATA}. 
To determine a {\bf non-overlapping count},
one reads the text from left to right.
Every time 
a match with the pattern is encountered, the matched word and all
characters to its left are removed before continuing the count.
For instance, there are 
2 non-overlapping occurrences of {\it  TATA}, 3 non-overlaping occurrences of {\it ATA},
and 1 non-overlapping occurrence of
{\it TA\#...TATA} in the text {\it ATATATATATA}.

\subsection{Automata, Probability, and Counting}
\label{sec:autprobcoun}

Important early work in pattern matching was done by Aho and Corasick
\cite{AhoCor75}, who constructed an automaton (now known as the
Aho-Corasick automaton) 
to search 
for a finite set of keywords in a text. Their work was focused on
bibliographic search and was therefore deterministic. The Aho-Corasick
automaton is an example of a deterministic finite automaton (DFA) that
we describe in detail in section~\ref{subsec:ahocor}.

We can imagine an Aho-Corasick automaton that recognizes the word {\it abba}.  Such an automaton would contain five states, numbered 1 (the empty string), 2 ($a$), 3 ($ab$), 4 ($abb$), and 5 ($abba$). A text is processed by the automaton one letter at a time, from left to right.  The automaton stays in state 1 until the presence
of an $a$ in the text triggers a transition to state 2. If the
next letter in the text is $b$, the automaton would then move to
state 3; otherwise the automaton would remain in state 2. (See figure
\ref{fig:abba} for a sketch of the transition rules for this
automaton.) If the automaton is in state 4 and the next letter
encountered is $a$, then the automaton would transition to state 5
which is associated with the detection of the keyword $abba$.

In automata used for pattern matching, prefixes of a word that are
also suffixes of the word play an important role in determining the
structure of the automaton. This structure is used in
mathematical techniques to count the number of words that contain or forbid certain patterns. The \textbf{autocorrelation of a string} and more generally
the \textbf{correlation between two strings} introduced by Guibas and Odlyzko quantifies this idea~\cite{GuiOdl81b}. The autocorrelation of a string $x$ is a string of 0's and 1's
--- of the same length as $x$ --- which gives information about the matches of $x$ with itself. The autocorrelation of $x$ is denoted $\Aut[x]$ and a 1 occurs at position $n$ in $\Aut[x]$ if and only if $x$ has a prefix of length $n$ which is also a suffix of $x$. For instance, if $x=abbab$ then $\Aut[x]=01001$. The {\bf autocorrelation polynomial} of a string $x$, denoted $\Aut[x;z]$, is the polynomial in the variable $z$ obtained by summing up all the monomials of the form $z^{n-1}$ for which the $n$-th character of $\Aut[x]$ is a 1. For instance, if $x=abbab$ then $\Aut[x;z]=z+z^4$. If $x$ is a string constructed with characters in an alphabet of size $s$ and $f(n)$ is the number of strings of length $n$ that do not have any occurrence of $x$ as a substring then
\[\sum_{n=0}^\infty\frac{f(n)}{z^n}=\frac{z\cdot\Aut[x;z]}{1+(z-s)\cdot\Aut[x;z]}.\]
See~\cite{GuiOdl81b} for generalizations of the above identity to consider more than just one forbidden strings.

Further discussion of automata theory is found in references \cite{HopUll79,Sip96,CroRyt02}. Combinatorial theory of pattern matching is discussed in references \cite{GouJac04,FlaSed06}. A thorough discussion of autocorrelation polynomials can be found in~\cite{Lot05}.

\subsection{Probabilistic counting}

The use of automata for pattern matching can be extended to consider
probabilities of occurrences of strings in random texts through the
use of probabilistic automata. We can think of a probabilistic
automaton
as a DFA that scans a random text
as the text is generated.
Transitions between different states of 
the
automaton occur according to probabilistic rules, which are determined
from the model which generates the random text (typically a Markov
chain).  
The probability that a word
occurs in a random text of a certain length can
be determined from the probability that 
the probabilistic automaton visits a specified set of states within a
certain number of steps.  This translates the pattern matching problem
into a problem about the behavior of a Markov chain. This
correspondence is helpful because the theory of Markov chains is a
well-established area 
of probability theory.

The occurrence of regular patterns in random strings
produced by Markov chains
(and more generally Hidden Markov chains) reduces
to problems regarding the behavior of a first-order homogeneous Markov
chain in the state space of a suitable DFA. This transformation of the
problem is
often called an \textbf{embedding technique}. As
we discuss in sections~\ref{sec:markov}, \ref{sec:app1}
and~\ref{sec:app2}, the embedding technique provides a unifying
theoretical framework for many different problems in probabilistic
pattern matching.

The \textbf{Markov chain embedding} technique 
usually corresponds to the embedding of
a random string into the states of an Aho-Corasick automaton.
In this framework, more complicated patterns 
(such as modular correlated patterns)
can be treated through the synchronization of
Aho-Corasick
automata 
associated with each possible combination of correlations. 
A \textbf{synchronized automaton} or \textbf{product automaton}
is a new automaton made up of multiple automata which simultaneously
process a single text. The technique of synchronization is discussed in detail in section~\ref{subsec:synch}.

Important early work in 
probabilistic pattern matching
was done by Li, who
studied the first occurrence of a reduced compound
pattern~\cite{Li80}.
A compound pattern is 
\textbf{reduced} if no word in the pattern is a substring of another
word in the pattern.
This work focused on a random text produced by a memoryless source.
Follow up work
by Gerber and Li 
studied
the
probability of occurrence of a reduced compound pattern in a random
string produced by a memoryless source~\cite{GerLi81}. Their approach
is based on martingale methods (of the type introduced 
in \cite{Li80})
and Markov chain embedding techniques which 
implicitly use
automata and synchronization. 
The martingale method developed
by Gerber and Li has been extended by Pozdnyakov and Kulldorff~\cite{PozKul06} in the setting of reduced compound patterns and memoryless sources but without the use of the Markov chain embedding technique.
 
Perhaps the most general
computational
treatment 
of pattern frequencies in random sequences was carried out by
Nicod\`eme, Salvy, and Flajolet \cite{NicSalFla02}. They considered
random strings produced either by a Bernoulli or a Markovian model and
focused on regular patterns which are of a non-degenerate form.  
A regular pattern is 
non-degenerate if the DFA that recognizes the pattern is irreducible
(i.e., from any state it is possible to visit any other state) and primitive (i.e., there is a minimal length $l$ such
that for any two pair of states in the automaton there exists a path
of length $l$ that connects the two states). 
Their analysis is based on automata theory and transfer matrix
methods.
They obtained the
generating function associated with the distribution of the number of
occurrences of a regular pattern in a random text.
From the generating functions they computed the mean and standard deviation of the Gaussian distribution associated with the number of occurrences of the pattern in sufficiently long sequences.  
When such an approach is applied to biological sequence analysis,
it allows the determination of $z$-scores associated with
different patterns,
and this allows researchers
to assess the significance of
matches.

Follow-up work by Nicod\`eme used automata theory and generating
functions as the basis of a symbolic package called
Regexpcount~\cite{Nic03}. This software can be used to study the
distribution of the number of occurrences of various regular
expressions in
Bernoulli or Markovian sources,
including simultaneous counts of different motifs. The software can
also perform searches for strings at a given {\bf edit distance} from a compound pattern and
compute the {\bf sooner-time} of a string,
given a random string with a prescribed prefix. The implementation of
the automata used
in this package relies on the concept of Marked automata
\cite{NicSalFla02}
and synchronization ideas. 
(We will not extensively discuss  Marked automata, but they can be
used as an alternative to synchronization.)

Another important reference in probabilistic pattern matching is the
book by Fu and Lou \cite{FuLou03}. This work compiles and extends
results of J. C. Fu and coauthors on the Markov chain embedding technique~\cite{FuKou94,FuCha02,FuCha03}. 
Although automata are not explicitly
used in these papers,
their embedding technique is effectively an implementation of the
Aho-Corasick automaton.  Their technique is applicable to the
occurrence or frequency of some compound, possibly modular patterns;
however, it 
cannot handle arbitrary regular patterns.  
Because their calculation technique typically requires a large number
of states, it has limited computational feasibility.

R\'egnier and coauthors have made important contributions to
probabilistic pattern matching
\cite{RegSzp98,Reg00b,RegLifMak00,RegDen04}. R\'egnier and
Szpan\-kowski studied overlap counting of a simple pattern;
their work considered
a text
generated by a first-order, stationary Markov chain \cite{RegSzp98}.
Their approach can be used to calculate generating functions
using a techinque which relies on
combinatorial relationships between certain languages (sets of words)
built from the pattern. 
They obtained relatively explicit forms for the generating functions,
in which the autocorrelation polynomial of the pattern being studied
appears naturally. As a result, they could extract the asymptotic
behavior of the coefficients that lead to central and large deviation
approximations for the distribution of the frequency statistic of the pattern.

In later work, R\'egnier generalized to $k$-th order stationary Markov sequences,
compound patterns, and either overlap or non-overlap counting
\cite{Reg00b}.  The paper gives insight into an aggregation procedure
of the words in a compound pattern that considerably simplifies the
complexity of the problem.  She defined \textbf{minimal languages}
associated with patterns, which contain no redundancies.
(This concept is distinct from the idea of a minimal automaton).
R\'egnier showed that the generating functions associated with the
minimal languages are determined by the generating functions
associated with some simpler auxiliary languages --- this allows an
important simplification of the calculations. The computation of
expectations, variances and correlations for the number of occurrences
of the different words in the compound pattern can be expressed
explicitly in terms of these generating functions. Her method is more
computationally efficient than some other
approaches that use automata to perform the same calculations, provided that the random string is produced by a stationary Markov source.

Two papers by R\'egnier and coauthors studied the over- and
underrepresentation of patterns.  R\'egnier,
Lifanov, and Makeev focused on compound patterns that are invariant
under the reverse-complement operation; they were studying the
counting of binding sites in double-stranded DNA \cite{RegLifMak00}.  This paper calculates $z$-scores to assess the over- or underrepresentation of patterns in
random sequences.  More recently, R\'egnier and Denise examined how
over- or underrepresentation of a pattern can depend on the over- or
underrepresentation of a second pattern (because information about the
frequency of the second pattern modifies the distribution of the first
pattern) \cite{RegDen04}.  In this paper, they studied the asymptotic
fraction of times that a single pattern is found in a random
string produced by a memoryless source or a stationary Markov source
of order $k$. The result is a large-deviation principle with an
explicit rate function and accompanying second-order local expansion.
The asymptotic expectation and standard deviation of a pattern conditioning on
the observed sequence of another pattern were determined in a
computable way.

Aston and Martin studied the probability that any of a set of compound
patterns is the first to be completed a certain number of times
\cite{AstMar05}. They studied binary strings produced by a Markovian
source. Their method is based on a Markov chain embedding technique
and allows the possibility that the count may be different for
different compound patterns.

Flajolet, Szpankowski and Vall{\'e}e studied the total number of
occurrences of a hidden pattern in a random text generated by a
memoryless source \cite{FlaSzpVal06}. A \textbf{hidden pattern}
appears in a string if all the characters in the pattern appear in
order in the string, although other arbitrary characters may appear
between the characters in the pattern. For example, the text {\it
  adenosine guanine} contains the hidden pattern {\it dog} because the
letters {\it d}, {\it o}, and {\it g} appear in order in the text.
In this paper, Flajolet and coauthors derive central limit theorems
for the number of occurrences of the hidden pattern, using a technique based on generating
function methods.
However, for what they call the fully constrained case
(i.e., when the gaps between letters in the hidden pattern are
constrained to be less than specified finite constants) they utilize
{\bf De Bruijn graphs} and transfer matrix methods to obtain more
refined results regarding the asymptotic distribution of the frequency
statistic. Recently, Bourdon and Vall\'ee have extended the analysis
of the asymptotic behavior of the expected value and variance for number of matches of a hidden pattern in a text generated by a dynamic source \cite{BouVal02}.  A \textbf{dynamic
  source} is a generalization of a Markovian source, where the
probability of a character may depend on all the preceding characters
\cite{Val01,CleFlaVal01}. 
Since Bernoulli and Markovian sources are special cases of dynamic
sources, use of dynamic source models is the most general theoretical
framework to study patterns in random strings.
Bourdon and Vall\'ee also showed that the frequency statistic
associated with a regular pattern in a random text produced by a
dynamical source is asymptotically Gaussian~\cite{BouVal06}.

\subsubsection{Forbidden patterns.}

Early work by Guibas and Odlyzko studied forbidden patterns
\cite{GuiOdl81b}. This paper addressed the probability that a reduced compound pattern does not appear in a random string produced by a memoryless source.  They
introduced the concept of autocorrelation polynomial to find the generating function associated with this probability. In probabilistic pattern matching, the \textbf{autocorrelation polynomial} of a string is also constructed from its autocorrelation (see section~\ref{sec:autprobcoun}) but taking into account the probabilities associated with the alphabet characters. See section 3 in~\cite{GuiOdl81b} for more details.

Gani and Irle also studied forbidden patterns~\cite{GanIrl99}. They determined the probability that a string of a given length does
not contain a 
type of compound pattern. 
The patterns they considered must be specified either by a completely
repetitive system or a system with a distinctive beginning (see their paper
for precise definitions).
Their approach is primarily computational and based on matrix
recursion methods. 
Their method is applicable to a memoryless source or a single
forbidden string in a text produced by a Markovian source.  In the
case of a Markovian source,
they constructed an automaton which is similar to the Aho-Corasick
automaton.

\subsubsection{Generalized words.}

Bender and Kochman studied the number of occurrences of generalized
words~\cite{BenKoc93}. They define a \textbf{generalized word} as a
set of strings of the same length. They focused on a memoryless source
and obtained central and local limit theorems for the joint
distribution of the number of occurrences of generalized words given
that a forbidden generalized word does not occur within the random
string; they were able to obtain explicit formulae only when there are
no forbidden generalized words. The use of de Bruijn automata and
transfer matrices is implicit in their argument.

\subsection{Distance between pattern occurrences}

These papers address the question of the separation between patterns
in the text; typically they are interested in computing the
probability that a pattern first occurs after $l$ characters of the
text or the probability that two patterns are separated by $m$ characters. 

The \textbf{sooner-time} of a pattern is the number of characters that precede the first occurrence of the pattern. Li calculated the expected value of the sooner-time of a reduced compound pattern \cite{Li80}.
This work focused on a random string produced by a memoryless source.
This approach is based on martingale techniques and also includes
calculation of the probability that any of the strings in the compound
pattern is first to occur.

Early work on sooner-times was motivated by the digestion of DNA by
restriction enzymes. In this experimental protocol, specific enzymes
recognize particular DNA sequences, called restriction sites; the
enzymes cut the DNA at the restriction sites. Typically the
restriction site can be described by a compound pattern, and one is
only interested in non-overlapping occurrences. This is justified
because enzymes cut the strand at the first position where a string in
the compound pattern is identified. Breen, Waterman, and Zhang found
the generating functions for this problem, assuming a random string
produced by a memoryless source \cite{BreWatZha85}. Their analysis is
based on renewal theory arguments \cite{Fel68} and autocorrelation
polynomials similar to those used in \cite{GuiOdl81b}. Biggins and
Cannings addressed the more general problem of Markovian sources
\cite{BigCan87}.

Robin and Daudin determined the exact distribution of (and generating
functions associated with) the distance between two consecutive
(possibly overlapping) occurrences of a reduced compound pattern
\cite{RobDau01}. They considered a random string produced by a 
first-order homogeneous Markov chain.  Their analysis is related to
autocorrelation polynomials; the technique is applied to analyze the
CHI-motif in the genome sequence of Haemophilus influenza.

Han and Hirano studied the distributions of sooner- and later-time for
two reduced patterns in a random string produced by a first-order
Markov chain \cite{HanHir03}. The \textbf{later-time} of two patterns is the number of characters that precede the completion of both patterns.
Their paper uses probabilistic arguments
to determine the generating functions associated with the sooner- and
later-time; their approach is related to the concept of
autocorrelation \cite{GuiOdl78,GuiOdl81a,GuiOdl81b}.  They also study
other statistics such as the distance between two successive
occurrences of the reduced patterns. Their argument can be adapted to
study the sooner-time of a reduced compound pattern.

Work by Park and Spouge studied the sooner-time and the distance between
occurrences for the more general case of a random text produced by an
irreducible, aperiodic stationary Markov chain \cite{ParSpo04}. This
approach used a Markov chain embedding technique (and implicitly
the Aho-Corasick automaton). They obtained in closed form the
generating function associated with the sooner-time and with the
statistic of distances between two consecutive occurrences of a
reduced compound pattern.

\subsection{Related techniques}

\subsubsection{Sequence alignment and seed sensitivity.}

Buhler, Keich, and Sun used techniques from automata theory and Markov
chains to determine optimal seeds for sequence alignment
\cite{BuhKeiSun03}.  Seeds are short
strings which are used as
starting points in sequence alignment algorithms to reduce the
computation time.
The approach of Buhler et al. allows the design of seeds that are
optimal (with respect a specified Markov model). They used the concept
of a similarity, which is used to quantify the matches between
sequences in an alignment. Their technique is based on a Markov chain embedding argument over the state space of an appropriate Aho-Corasick automaton.

Martin studied the distribution of the total number of successes (1s)
in success runs (sequences of 1s) longer than a predetermined length
in a binary sequence (sequence of 0s and 1s) produced by a Markov
source \cite{Mar05}. This work used a Markov chain embedding
technique.  It is applicable to the detection of tandem repeats in DNA
sequences: in this case a 1 corresponds to a match between two aligned
DNA sequences and a 0 to a mismatch. The distribution of the number of
successes is needed in the detection phase of Benson's
tandem-repeats-finder-algorithm \cite{Ben99} to validate candidate
sequences via hypothesis testing.

Kucherov, Noe, and Roytberg used automata theory to address the
general problem of determining seed sensitivity \cite{KucNoeRoy06}. In
this paper, Kucherov et al. permit the set of allowed seeds and target
alignments to be described by a DFA and allow the probabilistic model of the target alignments to be
described by a Hidden Markov model (rather than a finite-order Markov
chain). Their technique relies on a synchronization argument that
involves two DFAs and the HMM. They also define a new automaton to
specify the seed model that, according to simulation data, performs
2-3 orders of magnitude better than the
Aho-Corasick automaton.

\subsubsection{Random number generators.}

Work by Flajolet, Kirschenhofer, and Tichy studied the distribution of
substrings in binary strings \cite{FlaKirTic88}. Although the
motivation for this work is the performance of random number generators, the techniques used overlap with the techniques of pattern analysis in random strings.
They showed that almost all binary strings of length $n$ contain all
possible binary strings of length slightly less than $\log_2(n)$ a
nearly uniform number of times.  Their analysis is based on De Bruijn
graphs, auto-correlation polynomials, and generating functions.

\section{Languages, automata, and synchronization}
\label{sec:automata}

In this section, we introduce mathematical notation and definitions to
describe regular languages (section \ref{subsec:reglang}), automata
(section \ref{subsec:dfa}), and synchronization (section \ref{subsec:synch}). We
finalize with a discussion about Aho-Corasick automata (section \ref{subsec:ahocor}).
This section gives a self-contained presentation of the key
mathematical results and proofs for automata used in deterministic
pattern matching. 

\subsection{Main notation}

The {\it alphabet} $\A$ is a finite non-empty set; the elements in
$\A$ are characters used to construct strings. A {\it string} over
$\A$ is a finite sequence of characters in $\A$. We use lowercase
letters (such as $x$) to denote generic strings. The {\it length} of a
string $x$, denoted $|x|$, is the total number of characters (counting
all repetitions) in the string.  The {\it empty string}, denoted
$\epsilon$, is by definition the only string of length zero. We assume
that $\epsilon\notin\A$, that is, the alphabet does not contain the
empty string.  

The set $\A^*$ is defined to contain the empty string as well as all
strings formed with characters in $\A$.  A basic operation between two
strings is {\it concatenation}: if $x,y\in\A^*$ then $xy$ is defined
to be the string formed by concatenating $y$ after $x$. Since, by
definition, $x\epsilon=x$ and $\epsilon x=x$, in general
$|xy|=|x|+|y|$.

For $x\in\A^*$ and $1\le i\le j\le|x|$, $x[i..j]$ denotes the
substring of $x$ formed by all characters between and including the
$i$-th and $j$-th character of $x$. We write $x[i]$ as a shorthand for
$x[i..i]$. Note that for $x,y\in\A^*$, we write
$x\,y[i..j]$ to refer to the string formed by concatenating $x$ with
$y[i..j]$ as opposed to $(xy)[i..j]$ which refers to a substring of
$xy$.

For $x,y\in\A^*$ we write $x=...y$ to mean that there exists
$z\in\A^*$ (possibly empty) such that $x=zy$. In this case we say that
$y$ is a {\it suffix} of $x$. Similarly, we write $x=y...$ to mean
that there exists $z\in\A^*$ such that $x=yz$ and we say that $y$ is a
{\it prefix} of $x$.

\subsection{Regular Languages}
\label{subsec:reglang}

A {\it language} over $\A$ is any subset of $\A^*$; we typically use
the $\L$ to denote a language so $\L \subset \A^*$. We write $|\L|$ to refer to the {\it cardinality of $\L$}, i.e., the number of strings contained in $\L$. For example, $|\A|$ is the number of alphabet characters. This is not be confused with the length of a string: for $x\in\A^*$, $|x|$ refers to the length of $x$, however, $|\{x\}|=1$ regardless of the length of $x$ because $\{x\}$ is a language consisting of a single string.

Three standard operations, {\it union}, {\it concatenation} and {\it
  star}, are usually defined over languages. For
$\L_1,\L_2\subset\A^*$, the union $(\L_1\cup\L_2)$ corresponds to the
usual union of two sets, i.e., a string $x\in(\L_1\cup\L_2)$ if and
only if $x\in\L_1$ or $x\in\L_2$.  The concatenation language
$\L_1\L_2$ consists of all those strings of the form $xy$, with
$x\in\L_1$ and $y\in\L_2$.  Finally, $\L_1^*$ is the language formed
by the empty string and by any string that can be formed by
concatenating a finite number of strings in $\L_1$.  Mathematically,
$\L_1^*=\{\epsilon\}\cup\L_1\cup\L_1\L_1\cup\L_1\L_1\L_1\cup\ldots$

The class of {\it regular languages} is the smallest class of subsets
of $\A^*$ that contains all finite languages (i.e., languages
consisting of a finite number of strings) and that is closed under the
three standard operations.

\subsection{Deterministic Finite Automata}
\label{subsec:dfa}

A {\it deterministic finite automaton} (DFA) is a 5-tuple of the form
$G=(V,\A,f,q,T)$, where $V$ is a nonempty set, $\A$ is an alphabet,
$f:V\times\A\to V$ is a function, $q\in V$ and $T\subset V$. The terms
$V$, $f$, $q$ and $T$ are called, respectively, the {\it set of
  states}, {\it transition function}, {\it initial state} and {\it set
  of terminal states}.

In what follows $G=(V,\A,f,q,T)$ is a given DFA. $G$ can be
represented as a graph with vertex set $V$ where a directed edge
labeled with the character $\alpha$ goes from a vertex $u$ to a vertex
$v$ if and only if $f(u,\alpha)=v$. In particular, each vertex has
out-degree $|\A|$ and for all $u\in V$ and $\alpha\in\A$ there exists
a unique edge labeled with the character $\alpha$ that starts at $u$.
See figures \ref{fig:abba}, \ref{fig:baabba}, \ref{fig:aastarba} and \ref{fig:onebstara}
for examples of automata represented as directed labeled graphs.

The visual representation of $G$ facilitates the extension of the
transition function $f$ to the larger domain $V\times\A^*$ as follows.
For $x\in\A^*$ define the {\it path associated with $x$ in $G$ when
  starting at $u$} to be the sequence of states that are visited from
$u$ by following the edges in $G$ according to the labels appearing in
$x$ as they are read from left to right. In the special case that
$u=q$ (i.e., the path begins at the initial state), we refer to this path as
the {\it path associated with $x$ in $G$}. We define $f(u,x)$ to be
the state in $V$ where the path associated with $x$ ends when starting
at $u$. Note that $f(u,\epsilon)=u$. As a result, $f:V\times\A^*\to V$
satisfies the following fundamental property: for all $u\in V$ and
$x,y\in\A^*$,
\begin{equation}\label{ide:property f}
f(u,xy)=f(f(u,x),y).
\end{equation}
In other words, the path associated with the concatenation of two strings can be
determined by concatenating the paths associated with each string, provided that the end of the first path is used as the starting point of the second path.

For $u,v\in V$, we say that {\it $v$ is accessible from $u$} if there
exists $x\in\A^*$ such that $f(u,x)=v$.

The {\it language recognized by $G$} is defined as
\[L(G):=\{x\in\A^*:f(q,x)\in T\}.\]
In other words, $L(G)$
consists of all strings that can be formed by
concatenating from left to right the labels of the edges visited by
any path that starts at the initial state of $G$ and ends at some
terminal state.

In what follows we say that a language {$\L$ is recognized by $G$} provided that $\L=L(G)$. According to two
classical results in computer science, Kleene's theorem and the Rabin
and Scott theorem, the following holds \cite{HopUll79,Sip96}.

\begin{theorem}
  Let $\L\subset\A^*$. $\L$ is a regular language if and only if there
  exists a DFA $G$ such that $L(G)=\L$.
\end{theorem}

Consider two DFAs $G_1=(V_1,\A,f_1,q_1,T_1)$ and
$G_2=(V_2,\A,f_2,q_2,T_2)$. We say that {\it $G_1$ is isomorphic to
  $G_2$} (denoted $G_1\sim G_2$) provided that there is a bijection
$\Phi:V_1\to V_2$ such that $\Phi(q_1)=q_2$, $\Phi(T_1)=T_2$, and for
all $u,v\in V_1$ and $\alpha\in\A$, $f_1(u,\alpha)=v$ if and only if
$f_2(\Phi(u),\alpha)=\Phi(v)$. We can think of the function $\Phi$ informally as a relabeling of the states of $G_1$ that produces the states of $G_2$.  Using
(\ref{ide:property f}), one can see that $G_1 \sim G_2$ implies
that for all $u,v\in V$ and $x\in\A^*$
\[f_1(u,x)=v\Longleftrightarrow f_2(\Phi(u),x)=\Phi(v).\]
In particular, since $\Phi$ preserves initial states, the path
associated with $x$ in $G_1$ ends at $u$ if and only if the path
associated with $x$ in $G_2$ ends at $\Phi(u)$. Since $\Phi$ also
preserves terminal states, $G_1$ and $G_2$ recognize the
same language. Therefore, isomorphic automata recognize the same
regular languages.

\subsection{Synchronization}
\label{subsec:synch}

In what follows, for a given language $\L$, $\L^c$ denotes the {\it complement} of $\L$, i.e., $\L^c:=\{x\in\A^*:x\notin\L\}$.

Synchronization is an operation between two or more automata that can
be used to construct a new automaton that has useful properties, such
as recognizing multiple languages. To define synchronization, 
consider a finite sequence of regular languages $\L_i$,
$i=1,\ldots,m$, with $m\ge2$. For each $i$ let
$G_i=(V_i,\A,f_i,q_i,T_i)$ be a DFA that recognizes $\L_i$.

\begin{definition}
  The synchronized automaton associated with $G_1,\ldots,G_m$ is the
  automaton $G_1\times\cdots\times G_m=(V,\A,q,f,T)$ with
  $V:=V_1\times\cdots\times V_m$, $q:=(q_1,\ldots,q_m)$ and
  $T:=\{(u_1,\ldots,u_m)\in V:u_i\in T_i\hbox{ for at least one $i$}\}$. The transition function $f:V\times\A\to V$ is defined as
\begin{equation} f(\bfu,\alpha)
  :=(f_1(u_1,\alpha),\ldots,f_m(u_m,\alpha)),
\label{synchrule}
\end{equation}
for all $\bfu=(u_1,\ldots,u_m)\in V$ and $\alpha\in\A$.
To each $\bfu=(u_1,\ldots,u_m)\in V$ we associate the language
\[L(\bfu):=\left(\bigcap\limits_{i:u_i\in
      T_i}\L_i\right)\cap\left(\bigcup\limits_{i:u_i\notin
      T_i}\L_i\right)^c.\]
\end{definition}
Synchronized automata are also called {\it product automata}. We can
informally understand the idea of synchronization by imagining an
automaton which works by simultaneously operating the automata $G_1,
\ldots, G_m$: from the states $u_1,\ldots,u_m$ in the individual
automata, we feed each automaton the character $\alpha$. Then the
transitions of the synchronized automaton are determined by combining
all the transitions of the individual automata (which is what definition
(\ref{synchrule}) conveys). See figure \ref{fig:onebstara} for an example of a synchronized
automaton.

The key feature of synchronized automata is revealed by the following
result.

\begin{theorem}
\label{thm:synchronization}
If  $G_1\times\cdots\times G_m=(V,\A,q,f,T)$ then for all
$\bfu=(u_1,\ldots,u_m)\in V$ and $x\in\A^*$,
$f(\bfu,x)=(f_1(u_1,x),\ldots,f_m(u_m,x))$. In particular, for all
$x\in\A^*$, $x\in L(f(q,x))$.
\end{theorem}

\begin{proof}
  Fix $\bfu=(u_1,\ldots,u_m)\in V$. We show the first part by
  induction on the length of $x$. Since the case $|x|=0$ is trivial,
  it suffices to show that if the identity holds for all strings of
  length $n$ then it also holds for an $x\in\A^*$ of length $(n+1)$.
  Indeed, according to (\ref{ide:property f}), the inductive
  hypothesis and the definition of $f$, we have that
\begin{eqnarray*}
f(\bfu,x)&=&f(f(\bfu,x[1..n]),x[n+1]),\\
&=& f((f_1(u_1,x[1..n]),\ldots,f_m(u_m,x[1..n])),x[n+1]),\\
&=& (f_1(u_1,x),\ldots,f_m(u_m,x)),
\end{eqnarray*}
where we have used that $f_i(f_i(u_i,x[1..n]),x[n+1])=f_i(u_i,x)$ in the last identity. This proves the first part of the theorem.

For the second part, let $x\in\A^*$. According to the first part,
$f(q,x)=(f_1(q_1,x),\ldots,f_m(q_m,x))$. Since $G_i$ recognizes
$\L_i$, $x\in\L_i$ if and only if $f_i(q_i,x)\in T_i$. Consequently,
$x\in\cup_{i:f_i(q_i,x)\in T_i}\L_i$ and
$x\notin\cup_{i:f_i(q_i,x)\notin T_i}\L_i$. This completes the proof
of the theorem. $\Box$
\end{proof}

The first part of the theorem states that the path associated with $x$ in
the synchronized automaton is determined by the paths associated with
$x$ in each of the individual automata. A direct consequence of this is that $f(q,x)\in T$
if and only if there exists $i$ such that $f_i(q_i,x)\in T_i$.  In
other words, the synchronized automaton can reach a terminal state if
and only if one (or more) of the individual automata reaches a
terminal state. Since this is equivalent to having $x\in\L_i$, we see
that $G_1\times\cdots\times G_m$ recognizes the union language
$\cup_{i=1}^m\L_i$.

The second part of the theorem asserts that the state where the path associated with a string ends indicates all the languages $\L_1,\ldots,\L_m$ to which that string belongs to.
This permits to redefine the set of terminal states to recognize any language obtained via the intersections, unions and
complementations of the languages $\L_1,\ldots,\L_m$. For instance, if
we were to redefine $T$ as
\[\left\{\bfu\in V:
    L(\bfu)=\L_1^c\cap\left(\bigcup_{i=2}^m\L_i\right),\hbox{ or
    }L(\bfu)=\bigcap_{i=1}^m\L_i\right\}\]
then the resulting automaton would precisely recognize the language
\[\left(\L_1^c\cap\left(\bigcup_{i=2}^m\L_i\right)\right)
  \cup\left(\bigcap_{i=1}^m\L_i\right).\]
This feature of product automata is the key property used by computer
scientists to show that the class of regular languages is the same as
the class of languages recognized by DFAs (see
\cite{HopUll79,Sip96} for more details). In
pattern analysis in random sequences, this property is important for studying patterns that include but also exclude certain features.  

\subsection{Aho-Corasick automata}
\label{subsec:ahocor}

This class of automata was defined by Aho and
Corasick~\cite{AhoCor75} to detect all the occurrences
of a finite number of keywords in a general text.  Aho-Corasick
automata can be considered to be finite state machine implementations of the
Knuth-Morris-Pratt string searching algorithm~\cite{KnuMorPra77}.

\begin{figure}[t]\centering
 \includegraphics[height=6cm]{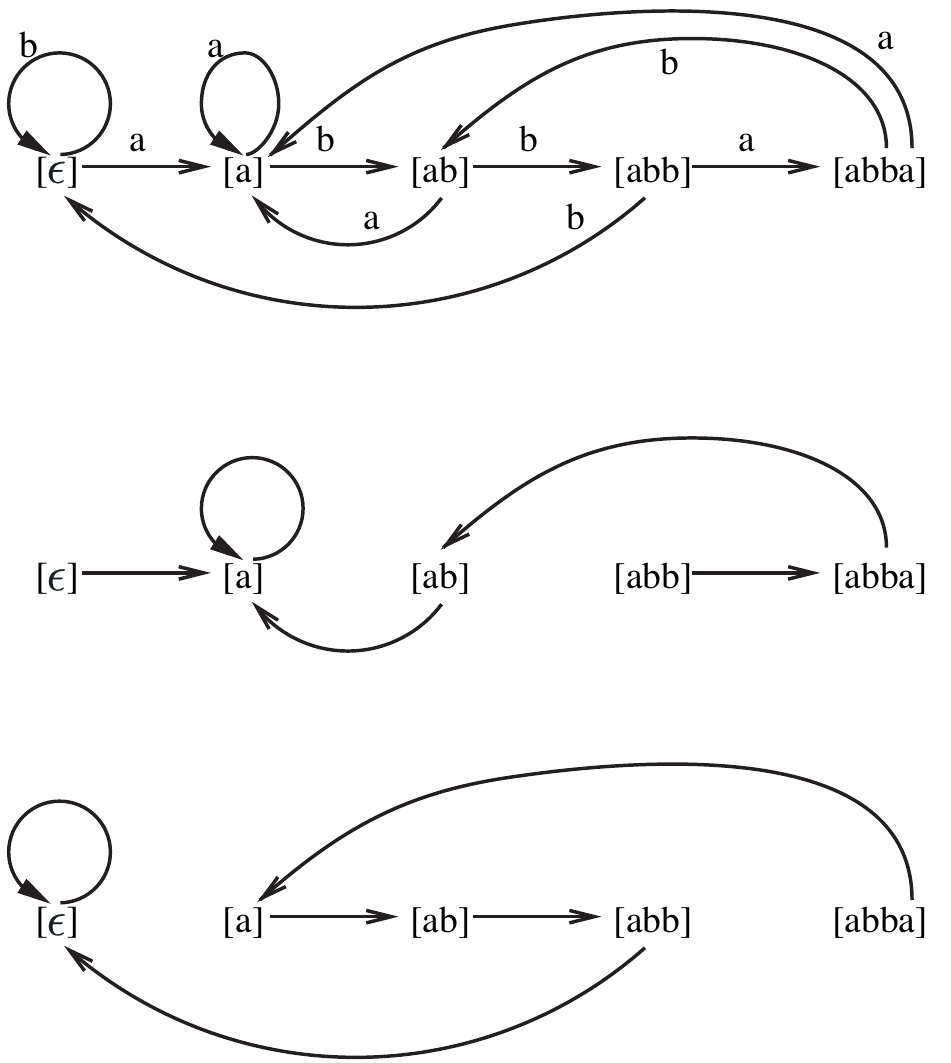}
 \caption[]{The Aho-Corasick automaton $AC(\{abba\})$ that recognizes the language $\{a,b\}^*\{abba\}$. Top, the full Aho-Corasick automaton that  finds all occurrences of $abba$ in a text constructed using the binary alphabet $\{a,b\}$. The initial state is the empty string (left). The terminal state is $abba$, which corresponds to detection of the string $abba$ in the text (right). Middle, the transitions that occur when the character $a$ occurs in the text. Bottom, the transitions that occur when the character $b$ occurs in the text.}
\label{fig:abba} 
\end{figure}

\begin{definition}
  Let $\W\subset\A^*$ be a finite non-empty set. The automaton
  $AC(\W)=(V,\A,q,f,T)$ is defined as follows. $V$ consists of the
  empty string as well as all prefixes of strings in $\W$,
  $q:=\epsilon$ and $T:=\W$. The transition function $f:V\times\A\to V$
  is defined such that for $u,v\in V$ and $\alpha\in\A$,
\[f(u,\alpha)=v\Longleftrightarrow \hbox{$v$ is the
    longest element in $V$ such that $u\alpha=...v$}.\]
\end{definition}
The main idea in the definition of the transition function $f$ is the
longest-prefix suffix rule: each state $u\in V$ contains information
about the longest prefix of a word in $\W$ that is at the same time a
suffix of a text so far scanned by the automaton. See figures \ref{fig:abba} and \ref{fig:baabba} respectively for a representation of $AC(\{abba\})$ and $AC(\{ba,abba\})$ as directed labeled graphs.

The technique used in~\cite{AhoCor75} to show the correctness of
Aho-Corasick automata relies on the concept of non-deterministic
finite automata. Here we present a new proof that is self-contained
and relies only on first principles. The following result can be
considered a rephrasing of Lemma 1 in \cite{AhoCor75}.

\begin{lemma}\label{lem:aho-corasick}
  For all $x\in\A^*$, $f(q,x)=u$ if and only if $u$ is the longest
  state in $V$ such that $x=...u$.
\end{lemma}

\begin{proof}
  We show the lemma by induction on the length of $x$. Since the case
  $|x|=0$ is trivial, it suffices to show that if $|x|=(n+1)$ and the
  lemma applies to all strings of length $n$ then it also applies for
  $x$. Let $u$ be the longest string in $V$ such that $x=...u$. Let
  $v=f(q,x[1..n])$ and $w=f(q,x)$. According to the inductive
  hypothesis, $v$ is the longest string in $V$ such that
  $x[1..n]=...v$. To prove the lemma it is enough to show that $u=w$.
  In order to do so we first show that
\begin{equation}\label{eq:|w|le|u|le|v|+1}
|w|\le|u|\le|v|+1\,.
\end{equation}
For this observe that according to (\ref{ide:property f}),
\begin{equation}\label{eq:edge x[n+1]}
f(v,x[n+1])=w.
\end{equation}
In particular, since $x=x[1..n] x[n+1]=...v x[n+1]$, it follows from
the above identity that $x=...w$. The defining property of $u$ implies
the first inequality in (\ref{eq:|w|le|u|le|v|+1}). To show the second
inequality, we proceed by contradiction.  Suppose, counterfactually,
that $|u|>|v|+1$. Since $x=...u$ and $x=...v x[n+1]$, there would be a
nonempty string $y$ such that $u=y v x[n+1]$. In particular, since
$u\in V$, $y v$ must be a prefix of a string in $\W$. Hence, $y v\in
V$. This is not possible because $x[1..n]=...y v$ and therefore $v$
could not be the longest element in $V$ with the property that
$x[1..n]=...v$. This contradicts the defining property of $v$ and
therefore the second inequality in (\ref{eq:|w|le|u|le|v|+1}) must be
true.

Finally, we show that $u=w$. Since $x=...u=...v x[n+1]$, the second
inequality in (\ref{eq:|w|le|u|le|v|+1}) implies that $v x[n+1]=...u$.
Using (\ref{eq:edge x[n+1]}), this implies that $|u|\le|w|$ and
therefore, according to the first inequality in
(\ref{eq:|w|le|u|le|v|+1}), $|u|=|w|$. Since $x=...u=...w$ then $u=w$.
This completes the proof of the lemma. $\Box$
\end{proof}

A direct consequence of the above lemma is that the Aho-Corasick
automaton $AC(\W)$ recognizes the language $\A^*\W$. However, its
terminal states satisfy an important property that is useful for
counting occurrences of patterns in random strings. The theorem
describing this property can be considered a rephrasing of Lemmas 2
and 3 in \cite{AhoCor75}.

\begin{theorem}\label{thm:aho-corasick}
  For $w\in\W$ define $T(w):=\{u\in\W: u=...w\}$. For all $w\in\W$ and
  $x\in\A^*$, $w$ occurs $m$ times as a substring of $x$ if and only
  if the path associated with $x$ in $AC(\W)$ visits the set $T(w)$
  exactly $m$ times.
\end{theorem}

\begin{proof}
  Suppose that $w$ occurs $m$ times as a substring of $x$ and that the
  path associated with $x$ in $AC(\W)$ visits $T(w)$ exactly $l$
  times. To prove the theorem it suffices to show that $m=l$.  Indeed,
  according to Lemma~\ref{lem:aho-corasick}, if for some $0\le
  i\le|x|$, $f(q,x[1..i])=u\in T(w)$ then $x[1..i]=...u=...w$.
  In particular, $m\ge l$. On the other hand, suppose that
  for some $1\le i\le|x|$, $x[1..i]=...w$. Let $u=f(q,x[1..i])$.
  According to the Lemma~\ref{lem:aho-corasick}, $u$ is the longest
  string in $V$ such that $x[1..i]=...u$. Since $w\in V$ and
  $x[1..i]=...w$, it follows that $|u|\ge|w|$. In particular, $u=...w$
  and therefore $u\in T(w)$.  This shows that $m\le l$ and hence
  $m=l$. This completes the proof of the theorem. $\Box$
\end{proof}

This theorem means that the Aho-Corasick automaton can be used to
count the number of occurrences of each of the keywords it searches for. In other
words, we can use the Aho-Corasick automaton to construct an automaton
that correctly matches any arbitrary set of strings $\W$. This
eliminates the need for the commonly used requirement in the analysis
of random strings that compound patterns be reduced.  A finite set of
strings $\W$ is said to be {\it reduced} provided that no string in
$\W$ is a substring of another string in $\W$. In this case,
$T(w)=\{w\}$ for each $w\in\W$ and therefore occurrences of $w$ in a
text are in one-to-one correspondence with the visits to state $w$ as
the automaton $AC(\W)$ processes the text. However, in order for this
last property to hold, it is enough that $\W$ is {\it suffix-reduced},
i.e., no string in $\W$ is a suffix of another string in $\W$. This
follows directly from Theorem~\ref{thm:aho-corasick} because for all
$w\in\W$, $T(w)=\{w\}$ precisely when $\W$ is suffix-reduced.

We finish this section with some remarks regarding the computational
complexity of Aho-Corasick automata. This type of automaton can be
implemented in time and space proportional to the sum of the lengths
of all words in $\W$. Furthermore, in the case of keyword sets with a
single string, Aho-Corasick automata turn out to be minimal: for all
$w\in\A^*$, $AC(\{w\})$ is the automaton with the smallest number of
vertices that recognizes the language $\A^*\{w\}$. 

We note that many algorithms other than Aho-Corasick can search for a
set of keywords. See \cite{HopUll79,CroRyt02} for more information. See \cite{Lot05} for an account of minimization algorithms that can be used to reduce the number of states of a given DFA.

\section{Markov chain embedding}
\label{sec:markov}

In this section, we extend the mathematical notation and definitions
introduced above to describe random walks on automata, a procedure
referred to as the Markov chain embedding. This procedure is the key
step required to move from deterministic to probabilistic pattern
matching, which is essential for the determination of the statistical
significance of genomic motif searches.  This section gives a
self-contained presentation of the key mathematical results and
proofs.

\subsection{Mathematical results}

As before, $\A$ is used to denote a generic alphabet. We introduce the concept of
a random text $X=(X_n)_{n\ge1}$, a sequence of $\A$-valued independent
and identically distributed random variables. The distribution of
$X_1$ in $\A$ is denoted as $\Prob(\cdot)$; in particular, for all
$n\ge1$ and $\alpha\in\A$, $\Prob(\alpha)$ corresponds to the
probability that $X_n=\alpha$.  We also define
\[\A^+:=\A^*\setminus\{\epsilon\}.\]
In other words, $\A^+$ is the set of all non-empty words formed by
concatenating characters in $\A$.

The following definition formalizes the notion of Markov chain
embedding as used in the literature by most authors.

\begin{definition}
  Let $G=(V,A,f,q,T)$ be a deterministic finite automaton. Define
  $V^G:=f(q,\A^+)$,i.e., $V^G$ is the set of all states in $u\in
  V$ for which there exists $x\in\A^+$ such that $f(q,x)=u$. The
  Markov chain embedding of $X$ in $G$ is the sequence of $V^G$-valued
  random variables $X^G:=(X^G_n)_{n\ge1}$ where
\[X^G_n:=f(q,X_1...X_n)\qquad(n\ge1).\]
\end{definition}

Recall that a sequence $Y=(Y_n)_{n\ge1}$ of $V^G$-valued random variables is said to be a first-order homogeneous Markov chain provided that for all $n\ge1$ and $u_1,\ldots,u_n,v\in V^G$,
\[P(Y_{n+1}=v\mid Y_n=u_n,\ldots,Y_1=u_1)=P(Y_{n+1}=v\mid Y_n=u_n),\]
and this last probability does not depend on $n$. The following theorem allows automatic computation of many statistics associated with patterns in random strings by connecting
the probabilistic calculations to the behavior of first-order
homogeneous Markov chains defined on the state space of an appropriate
automaton.

\begin{theorem}\label{thm:XG}
  If $X=(X_n)_{n\ge1}$ is a sequence of i.i.d. $\A$-valued random
  variables and $G=(V,A,f,q$, $T)$ is a deterministic finite automaton
  then $X^G$ is a first-order homogeneous Markov chain with initial
  distribution
\begin{equation}\label{thm:initial dist}
P(X^G_1=u)=\sum_{\alpha\in\A:f(q,\alpha)=u}\Prob(\alpha)\qquad(u\in V^G),
\end{equation}
and probability transitions 
\begin{equation}\label{thm:transition prob}
P(X^G_{n+1}=v\mid X^G_n=u) = \sum_{\alpha\in\A:f(u,\alpha)=v}
\Prob(\alpha)\qquad(u,v\in V^G). 
\end{equation}
\end{theorem}

\begin{proof}
  The proof of (\ref{thm:initial dist}) is direct. To show the Markov
  property observe that according to (\ref{ide:property f}),
  $X^G_{n+1}=f(X^G_n,X_{n+1})$. As a result, for all
  $u_1,\ldots,u_n,v\in V$ it applies that
\begin{eqnarray*}
&&P(X^G_1=u_1,\ldots,X^G_n=u_n,X^G_{n+1}=v)\\
&=& P(X^G_1=u_1,\ldots,X^G_n = u_n,f(u_n,X_{n+1})=v),\\
&=& P(X^G_1=u_1,\ldots,X^G_n=u_n)\cdot P(f(u_n,X_{n+1})=v),
\end{eqnarray*}
where for the second identity we have used that $X_{n+1}$ is
independent of $X_1,\ldots,X_n$.  This shows that $X^G$ is a
first-order Markov chain. Furthermore, since the distribution of
$X_{n+1}$ does not depend on $n$, it follows that the conditional
probability $P(X^G_{n+1}=v\mid X^G_n=u_n,\ldots,X^G_1=u_1)$ depends
only on $u_n$ and $v$ but not $n$. This shows that $X^G$ is
homogeneous. Therefore (\ref{thm:transition prob}) follows almost
immediately. This completes the proof.  $\Box$
\end{proof}

This theorem describes a random walk on the vertices of the automaton, where the
probability of a transition along an edge labeled with the character $\alpha$ is $\Prob(\alpha$). In other words, a
transition that occurs in the deterministic automaton in response to
reading character $\alpha$ occurs randomly with probability
$\Prob(\alpha)$. Therefore, the random walk can be represented by a
first-order Markov chain, where the transition probability depends
only on the current state and not on the preceeding states. A direct
consequence of this theorem is the following simple way to construct
the transition matrix of the Markov chain. To state the result we use {\it Iverson's brackets}: if $p$ is a statement then $[\![p]\!]=1$ provided that $p$ is a true statement, otherwise $[\![p]\!]=0$.

\begin{corollary}\label{cor:PG}
  If $G$ and $X$ are defined as in Theorem~\ref{thm:XG} then the
  probability transition matrix of $X^G$ in $V\times V$ is given by
  the formula
\begin{equation}\label{thm:prob trans matrix}
P^G=\sum_{\alpha\in\A}Prob(\alpha)\cdot G_\alpha,
\end{equation}
where $G_\alpha$ is the $V\times V$ matrix such that for all $u,v\in V$, $G_\alpha(u,v)=[\![f(u,\alpha)=v]\!]$.
\end{corollary}

In the above result, $G_\alpha$ corresponds to the incidence matrix of $G$ where only edges labeled with the character $\alpha$ are considered. See the middle and bottom part of figure~\ref{fig:abba} for a representation of $G_\alpha$ with $G=AC(\{abba\})$ and $\alpha=a\hbox{ or }b$.

\subsection{Prototype application of the Markov chain embedding}

Theorem~\ref{thm:XG} allows the calculation of the statistical
significance of matches of a regular pattern in a random string.  To
understand this application of the theorem, consider a random
text, i.e., a sequence $X=(X_n)_{n\ge0}$ of i.i.d.  random variables taking
values in some alphabet set $\A$. The patterns to be matched are
represented as a finite number of distinct regular languages
$\L_1,\ldots,\L_m$ in $\A^*$. We then define matches to each language
as
\[S^n_j:=\hbox{ number of substrings of $X_1...X_n$ that belong to $\L_j$}.\]
For each language (different $j$) we construct an automaton $G_j$ that
recognizes the regular language $\A^*\L_j$ and let $T_j$ denote the
set of terminal states of $G_j$. Define the synchronized automaton
constructed from the $G_j$, $G:=G_1\times\ldots\times G_m$ and let $T$
denote the set of terminal states of $G$. According to
Theorem~\ref{thm:synchronization}, there are $m_j$ (possibly
overlapping) substrings of $\L_j$ in $X_1...X_n$ provided that the
Markov chain $(X^G_i)_{i=1..n}$ visits the set of states $T(\L_j)$
exactly $m_j$ times. Therefore, if we define
\[T^n_j:=\hbox{ number of times that $(X^G_i)_{i=1..n}$ visits $T(\L_j)$}\]
then it follows that the vector of substring counts (the $S_j$) is
equal to the number of times the Markov chain visits the corresponding
terminal states:
\[(S^n_1,\ldots,S^n_m)=(T^n_1,\ldots,T^n_m).\]
In particular, the distribution of $(S^n_1,\ldots,S^n_m)$ can be
completely studied in terms of the distribution of
$(T^n_1,\ldots,T^n_j)$, to which we can apply the theory of Markov chains. 

Several refinements of the above method are possible for different
tasks. For instance, if we are interested in forbidden patterns, the
probability that no substring of $X_1...X_n$ belongs to
$\cup_{j=1}^m\L_j$ corresponds to the probability that
$(T^n_1,\ldots,T^n_j)=(0,...,0)$. In addition, the over- or
underrepresentation of patterns described by the languages $\L_1$ and $\L_2$ given the vector of counts for languages $\L_3,\ldots,\L_{m}$ could be studied
in terms of the joint distribution of $(T^n_1,\ldots,T^n_m)$ and the
marginal distribution of $(T^n_3,\ldots,T^n_m)$. Finally, the
aggregated number of occurrences of strings in $\cup_{j=1}^m\L_j$ as
substrings of $X_1...X_m$ corresponds to the total number of visits
that $(X^G_i)_{i=1..n}$ makes to $T$.

The particular form of the product automaton $G_1\times\cdots\times
G_m$ we have been using is sometimes not computationally efficient. Indeed, in many situations the product automaton has a computationally intractable number of states.  The key mathematical property of this automaton is that its states are associated with the detection or non-detection of each of the languages $\L_1,\ldots,\L_m$.  This allows one to determine the distribution of $(S^n_1,\ldots,S^n_m)$ in terms of the
Markov chain $X^G$. However, this property is not exclusive to product
automata. Other authors have proposed automata with similar
characteristics called {\it Marked automata} and that can
be used in the context of regular languages and random strings modeled by Markov sources \cite{NicSalFla02,Nic03}. For a related discussion see~\cite{Lla07} where a synchronization argument is used to construct the smallest state space size automaton required for analyzing the number of matches with a regular pattern in a random string generated by a Markov source.

\section{Application to a compound pattern}
\label{sec:app1}

This section considers a prototype example for studying the sooner-time and frequency statistic of a possibly non-reduced compound pattern in random strings produced by memoryless source. In this context, any string in a compound pattern counts as a match. Potential applications of this apparatus include the study of RNA motifs, in which the compound pattern might include a degenerate base (e.g.,  the
symbol $R$ stands for either of the two purines, $A$ and $G$, so the
sequence $CCRU$ represents the compound pattern $\{CCAU,CCGU\}$), or
by base pairing (e.g., the sequence $1GAAA1'$ --- with $A':=U$, $C':=G$, $G':=C\hbox{ or }U$ and $U':=A\hbox{ or }G$ --- allows the first and last nucleotide to pair with each other, the compound pattern is $\{AGAAAU,\,CGAAAG,\,GGAAAC,\,GGAAAU,\,UGAAAA,$ $UGAAAG\,\}$).

We will use two patterns on a binary alphabet to illustrate the main
principles. Consider the alphabet $\A=\{a,b\}$ and let
$X=(X_n)_{n\ge1}$ be a sequence of i.i.d. $\A$-valued random
variables with initial distribution $P(X_1=a)=p$ and $P(X_1=b)=q$,
with $p\cdot q>0$ and $p+q=1$.  In this example we study the
occurrences of the patterns $ba$ and $abba$ in $X$.

\begin{figure}[t]\centering
  \includegraphics[height=8cm]{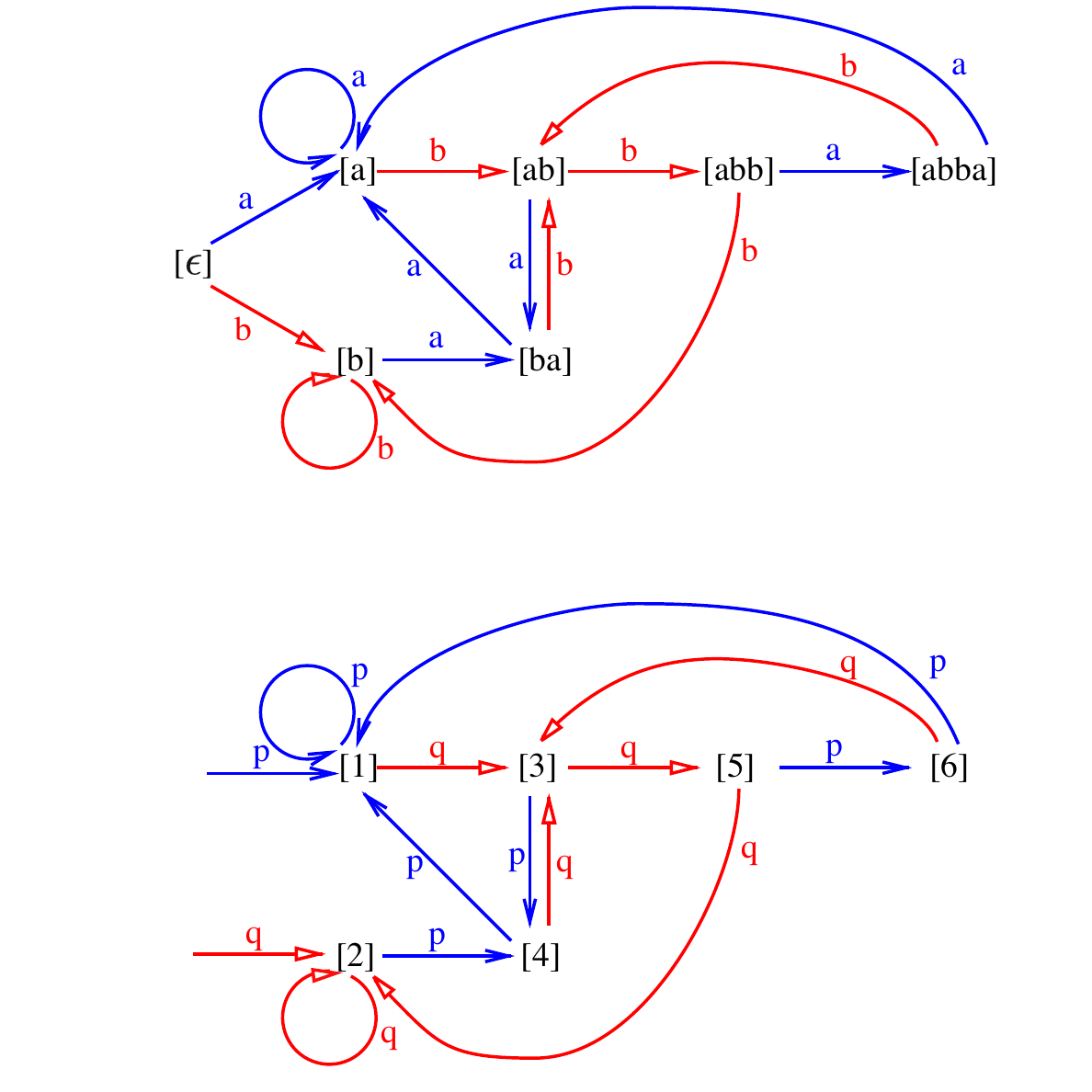}
 \caption[]{The automaton that recognizes the non-reduced compound pattern
$\{ba,abba\}$ in a text constructed using the binary alphabet $\{a,b\}$. Top, the
Aho-Corasick automaton $AC(\{abba,ba\})$ which detects all
occurrences of $ba$ and $abba$ in a binary text. The initial state is the empty string (left), and the terminal states are $ba$ and $abba$ (right).
Bottom, representation of the first-order homogeneous Markov chain associated with a random text embedded in the automaton on top. The text is produced by a memoryless source where the character $a$ occurs with probability $p$ and the character $b$ occurs with probability $q$. The Markov chain starts at state $1$ with probability $p$ and at state $2$ with probability $q$. The probability that $abba$ occurs in a random string of length $n$ is equivalent to the probability that the Markov chain visits state 6 within $(n-1)$ steps. Similarly, the probability that $ba$ occurs in a random string is equivalent to the probability that the Markov chain visits states 4 or 6.
}
\label{fig:baabba} 
\end{figure}

For the rest of this section, $G$ denotes the Aho-Corasick automaton
$AC(\{ba,abba\})$. A visual representation of $G$ is given in figure
\ref{fig:baabba}.

According to Theorem~\ref{thm:XG}, $X^G$ is a first-order homogeneous
Markov chain with states $a$, $b$, $ab$, $ba$, $abb$, $abba$ which we
label respectively as $1$, $2$, $3$, $4$, $5$, $6$. From
(\ref{thm:initial dist}), (\ref{thm:transition prob}) and
(\ref{thm:prob trans matrix}) it follows that $X^G$ has an initial
distribution given by the vector
\[\mu:=[\begin{array}{cccccc}
p & q & 0 & 0 & 0 & 0
\end{array}]\]
and the probability transition matrix
\[P:=p\cdot\!\left[\begin{array}{cccccc}
1 & 0 & 0 & 0 & 0 & 0\\
0 & 0 & 0 & 1 & 0 & 0\\
0 & 0 & 0 & 1 & 0 & 0\\
1 & 0 & 0 & 0 & 0 & 0\\
0 & 0 & 0 & 0 & 0 & 1\\
1 & 0 & 0 & 0 & 0 & 0
\end{array}\right]
+
q\cdot\!\left[\begin{array}{cccccc}
0 & 0 & 1 & 0 & 0 & 0\\
0 & 1 & 0 & 0 & 0 & 0\\
0 & 0 & 0 & 0 & 1 & 0\\
0 & 0 & 1 & 0 & 0 & 0\\
0 & 1 & 0 & 0 & 0 & 0\\
0 & 0 & 1 & 0 & 0 & 0
\end{array}\right]
=
\left[\begin{array}{cccccc}
p & 0 & q & 0 & 0 & 0\\
0 & q & 0 & p & 0 & 0\\
0 & 0 & 0 & p & q & 0\\
p & 0 & q & 0 & 0 & 0\\
0 & q & 0 & 0 & 0 & p\\
p & 0 & q & 0 & 0 & 0
\end{array}\right].\]
A visual representation of $X^G$ is displayed in figure \ref{fig:baabba}.

According to Lemma~\ref{lem:aho-corasick}, each occurrence of $abba$
in $X^G$ corresponds to a visit to state $6$. On the other hand, each
occurrence of $ba$ which does not contribute to an occurrence of
$abba$ corresponds to a visit to state $4$. In particular, all
occurrences of $ba$ in $X^G$ correspond to visits to states $4$ and
$6$.

\subsection{Sooner-time distribution of two non-reduced patterns}
\label{sec:sooner-time non red}

Broadly speaking, the sooner-time of a pattern corresponds to the position of the first occurrence of the pattern in a random text. Potential applications of our apparatus include the analysis of the occurrence of any one of a set of completely different RNA patterns that can catalyze the same reaction, such as RNA self-cleavage \cite{tang00}.

Define
\[T:=\hbox{  sooner-time distribution of $ba$ or $abba$},\]
i.e., $T$ is the smallest $n$ such that $X_1...X_n=...ba$ or
$X_1...X_n=...abba$. To study the distribution of $T$ consider the
matrix and vectors
\[Q:=
\left[\begin{array}{cccc}
p & 0 & q & 0\\
0 & q & 0 & 0\\
0 & 0 & 0 & q\\
0 & q & 0 & 0\\
\end{array}\right]
\quad;\quad
\nu:=[\begin{array}{cccc}
p & q & 0 & 0
\end{array}]
\quad;\quad
u:=\left[\begin{array}{c}
0\\
p\\
p\\
0
\end{array}\right]
\quad;\quad
v:=\left[\begin{array}{c}
0\\
0\\
0\\
p
\end{array}\right].\]
The matrix $Q$ is obtained by removing the fourth and sixth rows and columns of the probability transition matrix $P$, i.e., the entries associated with the patterns $ba$ and $abba$. The vector $\nu$ corresponds to the vector $\mu$ with the fourth and sixth columns removed. The vectors $u$ and $v$ correspond to the fourth and sixth columns of $P$ with the fourth and sixth rows removed.

Each entry in any power of $Q$ is an aggregate probability, the
probability that $X^G$ follows certain paths that avoid any edge that
is incident to states $4$ or $6$. The entry in row $r$
and column $c$ of $Q^n$ corresponds to the probability that $X^G_n=c$
and $(X^G_i)_{i=1,...,n}$ does not visit states $4$ and $6$, given
that $X^G_1=r$. Since the entry in row $r$ of $(u+v)$ corresponds to the probability that $X^G_{n+1}=4\hbox{ or }6$ given that $X^G_n=r$, it follows that
\begin{equation}\label{ide:P(T=n)}
\Prob[T=n]=\nu\cdot Q^{n-2}\cdot (u+v)\qquad(n\ge2).
\end{equation}

The above expression can be rewritten in terms of the
generating function of $T$. Let $\Id_4$ be the $(4\times 4)$
identity matrix. Since
\begin{equation}\label{ide:matrix series}
\sum_{n=0}^\infty z^n\cdot Q^n=(\Id_4-z\cdot  Q)^{-1}\qquad(|z|<1),
\end{equation} 
it follows from (\ref{ide:P(T=n)}) that
\[\sum_{n=2}^\infty \Prob[T=n]\,z^n=z^2\cdot\nu\cdot(\Id_4-z\cdot Q)^{-1}\cdot (u+v).\]

Using the cofactor formula to invert the matrix
on the right-hand side above, we obtain for this example that
\begin{equation}\label{eq:f(z)}
\sum_{n=2}^\infty\Prob[T=n]\,z^n =
  \frac{pqz^2}{(1-pz)(1-qz)}.
\end{equation}

Before continuing we introduce some standard notation~\cite{Wil94,FlaSed06}. In what follows, wherever $f(z)$ is a power series in the variable $z$, i.e., $f(z)=\sum_{n=0}^\infty f_nz^n$ with $f_0,f_1,\ldots$ complex numbers, the coefficient of $z^n$ of $f(z)$ is denoted $[z^n]f(z)$. Specifically, $[z^n]f(z):=f_n$. For instance, using a geometric series argument, it follows that
\[\frac{1}{\alpha-\beta z}=\sum_{n=1}^\infty\frac{\beta^n}{\alpha^{n+1}}z^n\qquad(\alpha,\beta\ne0;|z|<|\alpha/\beta|).\]
In particular, $[z^n]1/(\alpha-\beta z)=\beta^n\alpha^{-(n+1)}$. Via successive differentiation of both sides above with respect to the variable $z$, one obtains for all integer $m\ge1$ the following well-known formula~\cite{Wil94}
\begin{equation}\label{ide:coef (1-x)mm} 
[z^n]\frac{1}{(\alpha-\beta
  z)^m}=\frac{\beta^m\alpha^{-(n+m)}}{(m-1)!}
\prod_{j=1}^{m-1}(n+j)\qquad(\alpha,\beta\ne0;n\ge0). 
\end{equation}

To obtain an explicit formula for $\Prob[T=n]$ we use the
partial fraction decomposition of the right-hand side of equation
(\ref{eq:f(z)}) and then (\ref{ide:coef (1-x)mm}) to extract the coefficient of $z^n$ in each of the terms of the
decomposition~\cite{Wil94,FlaSed06}. For instance, if $p\ne q$ then the partial fraction decomposition of the right-hand side in~(\ref{eq:f(z)}) leads to the identity
\[\sum_{n=2}^\infty \Prob[T=n]\,z^n=\frac{q}{(1-pz)(p-q)}+\frac{p}{(1-qz)(q-p)}+1.\]
As a result, using (\ref{ide:coef (1-x)mm}) to identify the coefficient of $z^n$ in each of the terms on the right-hand side above, it follows that
\[\Prob[T=n]=\frac{q\cdot p^n-p\cdot q^n}{p-q}\qquad(n\ge2;p\ne q).\]
On the other hand, if $p=q$ then
\[\sum_{n=2}^\infty\Prob[T=n]\,z^n=\frac{4}{(2-z)^2}-\frac{4}{2-z}+1,\]
and therefore
\[\Prob[T=n]=\frac{n-1}{2^n}\qquad(n\ge2;p=q).\]

To study the probability of which of the patterns $ba$ or $abba$ is
the first to be observed, notice that by definition of $T$, $X_1...X_T=...ba$ or $X_1...X_T=...abba$. To determine the probability
that $ba$ is observed before $abba$, or the probability that $ba$ and
$abba$ are observed simultaneously for the first time we use that
\begin{eqnarray*}
P(T=n,X_1...X_T\ne...abba)&=&\nu\cdot Q^{n-2}\cdot u\qquad(n\ge2),\\
P(T=n,X_1...X_T=...abba)&=&\nu\cdot Q^{n-2}\cdot v\qquad(n\ge2).
\end{eqnarray*}
Since $\det(\Id_4-Q)=(1-p)\cdot(1-q)$, it follows from (\ref{ide:matrix series}) that
\begin{eqnarray*}
P(X_1...X_T\ne...abba)&=&\nu\cdot(\Id_4-Q)^{-1}\cdot u,\\
\nonumber
P(X_1...X_T=...abba)&=&\nu\cdot(\Id_4-Q)^{-1}\cdot v.
\end{eqnarray*}
Using symbolic algebra software to evaluate the right-hand side of the above identities
we find that the probability that $ba$ is observed before pattern
$abba$ is $(1-p^2q)$. The probability that $ba$ and $abba$ are
observed simultaneously for the first time is therefore $p^2q$.

\begin{table}[h]
\caption{Joint distribution for the frequency statistics $S_6^1$ and $S_6^2$ as defined in (\ref{ide:Sn1}) and (\ref{ide:Sn2}), respectively. Since the probabilities in the third column add up to one, no other combination of $m_1$ and $m_2$ is possible for a random binary string of length $6$. The probabilities in the third column can be computed via matrix multiplication using identity (\ref{ide:prob S_n12=m1m2}), or by extracting the coefficient of $x^6y_1^{m_1}y_2^{m_2}$ of the generating function $F(x,y_1,y_2)$ in (\ref{ide:gen fct Sn12}) or (\ref{ide:explicit gen fct Sn12}).}
\centering
\label{table:joint dist Sn12}
\begin{tabular}{lll}
\hline\noalign{\smallskip}
$m_1$ & $m_2$ & Probability that $(S_6^1,S_6^2)=(m_1,m_2)$\\[3pt]
 0 & 0 & $p^3q^3+p^6+p^5q+p^4q^2+p^2q^4+pq^5+q^6$\\[2pt]
 1 & 0 & $7p^3q^3+7p^2q^4+5pq^5+5p^4q^2+5p^5q$\\[2pt]
 2 & 0 & $5p^3q^3+4p^2q^4+6p^4q^2$\\[2pt]
 3 & 0 & $p^3q^3$\\[2pt]
 1 & 1 & $2p^3q^3+p^2q^4+3p^4q^2$\\[2pt]
 2 & 1 & $2p^2q^4+4p^3q^3$\\
\noalign{\smallskip}\hline
\end{tabular}
\end{table}

\begin{table}[h]
\caption{The joint distribution of $S_6^1$ and $S_6^2$ displayed in Table \ref{table:joint dist Sn12} permits the calculation of the conditional distribution of $S_6^1$ given $S_6^2=0$. In particular, the expected value and variance of the number of occurrences of $ba$ as a substring of $X_1...X_6$ can be reassessed when the pattern $abba$ is known to not occur as a substring of the random string.}
\centering
\label{table:conditional dist Sn12}
\begin{tabular}{ll}
\hline\noalign{\smallskip}
$m_1$ & Probability that $S_6^1=m_1$ given that $S_6^2=0$\\[3pt]
 1 & $\frac{5p-18p^2+29p^3-24p^4+13p^5-5p^6}{1-3p^2+6p^3-3p^4}$\\[2pt]
 2 & $\frac{5p^6-13p^5+15p^4-11p^3+4p^2}{1-3p^2+6p^3-3p^4}$\\[2pt]
 3 & $\frac{p^3-3p^4+3p^5-p^6}{1-3p^2+6p^3-3p^4}$\\
\noalign{\smallskip}\hline
\end{tabular}
\end{table}

\subsection{Frequency statistics of two non-reduced patterns}
\label{sec:freq non-reduced}

In this section we study the joint distribution of the number of
occurrences of the patterns $ba$ and $abba$ in $X_1...X_n$.
Observe that $\{ba,abba\}$ is not a reduced set of patterns because
$ba$ is a suffix of $abba$. One might be interested in non-reduced
patterns, for example, in studying the combinatorial function of
miRNA seed sequences, which may or may not act together
to regulate gene expression during translation \cite{lewi05}. In order
to test whether the occurrence of two seeds is correlated, we would
need to first calculate the null distribution if there were no
functional relationship: even in the absence of biological effects,
the probability of observing one seed might affect the probability of
observing the other (for example, if one seed were to overlap the
other). For this example, again demonstrated on the two-letter
alphabet, consider the random variables
\begin{eqnarray}
\label{ide:Sn1}
S_n^1&:=&\hbox{number of times that $ba$ occurs as a substring of $X_1...X_n$},\\
\label{ide:Sn2}
S_n^2&:=&\hbox{number of times that $abba$ occurs as a substring of $X_1...X_n$}.
\end{eqnarray}
The argument to be presented here could also be used to study the distribution of $(S_n^1-S_n^2,S_n^2)$, where $(S_n^1-S_n^2)$ corresponds to the number of times that $ba$ appears as a substring of $X_1...X_n$ but without contributing to an occurrence of $abba$ as a substring.

The notation introduced in section~\ref{sec:sooner-time non red} will be extended to consider power series in several variables. For instance, if $g(x,y)=\sum_{n,m=0}^\infty g_{n,m}x^ny^m$, with $(g_{n,m})_{n,m\ge0}$ an array of complex numbers, we define $[x^ny^m]g(x,y):=g_{n,m}$.

To study the joint distribution of $S_n^1$ and $S_n^2$ we use a
transfer matrix method~\cite{FlaSed06,GouJac04}.  Consider the matrix with polynomial
entries and the vector
\[P_{y_1,y_2}:=
\left[\begin{array}{cccccccccccccccc}
p & & & 0 & & & q & & & 0 & & & 0 & & & 0\\
0 & & & q & & & 0 & & & py_1 & & & 0 & & & 0\\
0 & & & 0 & & & 0 & & & py_1 & & & q & & & 0\\
p & & & 0 & & & q & & & 0 & & & 0 & & & 0\\
0 & & & q & & & 0 & & & 0 & & & 0 & & & py_1y_2\\
p & & & 0 & & & q & & & 0 & & & 0 & & & 0
\end{array}\right]
\quad;\quad
\delta:=\left[\begin{array}{c}
1\\
1\\
1\\
1\\
1\\
1\end{array}\right].\]
The matrix $P_{y_1,y_2}$ is obtained by multiplying the fourth and
sixth column of $P$ by $y_1$, and the sixth column of $P$ by $y_2$.
Observe that if $y_1=1$ and $y_2=1$ then the entry in row $r$ and
column $c$ of $P^n_{y_1,y_2}$ corresponds to the probability that
$X^G_n=c$ given that $X^G_1=r$. This is because the entries in $P^n$
correspond to the aggregate probability of all possible paths of
length $n$ that start at state $r$ and end at state $c$. The entry in row $r$ and column $c$ of $P^n_{y_1,y_2}$ is a polynomial in the
variables $y_1$ and $y_2$. The coefficient of
$y_1^{m_1}y_2^{m_2}$ is the aggregate probability of all paths of
length $n$ that start at $r$, end at $c$, and visit $m_1$ times the
set of states $\{4,6\}$ and $m_2$ times the set $\{6\}$. As a result,
for $n\ge1$ and $m_1,m_2\ge0$ one finds that
\begin{equation}\label{ide:prob S_n12=m1m2}
\Prob[S_n^1=m_1,S_n^2=m_2]=[y_1^{m_1}y_2^{m_2}](\mu\cdot
P^{n-1}_{y_1,y_2}\cdot\delta). 
\end{equation}
This suffices to determine the joint distribution of $S_n^1$ and
$S_n^2$ for small values of $n$. See tables~\ref{table:joint dist Sn12} and~\ref{table:conditional dist Sn12} for specific computations in the case of $n=6$.

Define
\[F(x,y_1,y_2):=\sum_{n=1}^\infty \sum_{m_1=0}^\infty
  \sum_{m_2=0}^{m_1} 
  \Prob[S_n^1=m_1,S_n^2=m_2]\,x^ny_1^{m_1}y_2^{m_2}.\]
In terms of generating functions, identity (\ref{ide:prob S_n12=m1m2}) figure
\begin{equation}\label{ide:gen fct Sn12}
F(x,y_1,y_2)=x\cdot\mu\cdot(\Id_6-x\cdot P_{y_1,y_2})^{-1}\cdot\delta,
\end{equation}
where $\Id_6$ is the $(6\times6)$ identity matrix. The matrix on the right-hand side above can be determined in closed form using symbolic algebra software. By doing so one derives that
\begin{eqnarray}
&&F(x,y_1,y_2)\nonumber\\
\label{ide:explicit gen fct Sn12}
&=&\frac{pq^3y_1(1-y_2)x^4+pq(y_1-1)x^2+x}{pq^3y_1(y_2-1)x^4+pq^2y_1(1-y_2)x^3+pq(1-y_1)x^2-x+1}.
\end{eqnarray}
According to the definition of $F(x,y_1,y_2)$, the coefficient
of $x^ny_1^{m_1}y_2^{m_2}$ on the right-hand side above corresponds to the probability that
$(S_n^1,S_n^2)=(m_1,m_2)$. For small values of $n$ this allows a
direct calculation of the joint distribution of $S_n^1$ and $S_n^2$ by determining the Taylor coefficients of $F(x,y_1,y_2)$ about $(x,y_1,y_2)=(0,0,0)$.

For large values of $n$ an asymptotic analysis of the joint distribution of
$S_n^1$ and $S_n^2$ is more appropriate. This follows in the general
context of linear (also called additive) functionals of Markov chains that we briefly describe next. For this consider an integer $d\ge1$. In what follows, $d$-dimensional vectors are thought of as column vectors. For a $d$-dimensional vector $c$, we write $c'$ to refer to the transpose of $c$. Consider a vector-valued transformation $f=(f_1,\ldots,f_d)'$, where each entry $f_i:V^G\to\RR$ is a given function. We are interested in the asymptotic behavior of the random variables
\[S_n^f:=\sum_{i=1}^n f(X^G_i).\]
Define the $d$-dimensional vector and $(d\times d)$ matrix
\begin{eqnarray}
\label{ide:CLT mu}
\mu&=&\lim_{n\to\infty}\frac{1}{n}\left[\begin{array}{c}
\E(S_n^{f_1})\\
\vdots\\
\E(S_n^{f_d})
\end{array}\right],\\
\label{ide:CLT Sigma}
\Sigma&=&\lim_{n\to\infty}\frac{1}{n}\left[\begin{array}{ccc}
\Var(S_n^{f_1}) & \ldots & \Cov(S_n^{f_1},S_n^{f_d})\\
\vdots & \ddots & \vdots\\
\Cov(S_n^{f_d},S_n^{f_1}) & \ldots & \Var(S_n^{f_d})
\end{array}\right].
\end{eqnarray}
Whenever $X^G$ is an aperiodic and irreducible first-order homogeneous
Markov chain in a finite state space, the entries
in $\mu$ and $\Sigma$ above are finite and do not depend on the initial
distribution of $X^G$~\cite{Che99,Jon04}. In particular, $\Sigma$ is a semi-positive
definite matrix, i.e., $c'\cdot\Sigma\cdot c\ge0$ for all $d$-dimensional vector $c$. 
Furthermore, the aperiodicity and irreducibility of $X^G$ implies
that $(S_n^f-n\mu)/\sqrt{n}$ converges in distribution to a centered
$d$-dimensional normal distribution with variance-covariance matrix
$\Sigma$. (This follows from the Cram\'er-Wold device~\cite{Sha03} and the general results in~\cite{Che99,Jon04}.) This means that for each pair of real numbers $a\le b$ and $d$-dimensional vector $c$ such that $c'\cdot\Sigma\cdot c>0$,
\begin{eqnarray*}\lim_{n\to\infty}P\left(a\le c'\cdot\frac{S_n^f-n\mu}{\sqrt{n}}\le b\right)&=&\frac{1}{\sqrt{2\pi(c'\cdot \Sigma\cdot c)}}\\
&&\hphantom{1234}\cdot\int_a^b\exp\left\{-\frac{x^2}{2(c'\cdot\Sigma\cdot c)}\right\}dx.
\end{eqnarray*} 
In addition, if $\det\Sigma>0$ then
\begin{eqnarray*}
\lim_{n\to\infty}P\left(\frac{S_n^f-n\mu}{\sqrt{n}}\in\Theta\right)=\frac{1}{(2\pi\cdot\det\Sigma)^{d/2}}\cdot\int_\Theta \exp\left\{-\frac{x'\cdot\Sigma^{-1}\cdot x}{2}\right\}dx,
\end{eqnarray*}
for all measurable sets $\Theta\subset\RR^d$ whose boundary $\partial\Theta$ is of Lebesgue measure zero.

In the context of the frequency statistics $(S_n^1,S_n^2)$ consider the function $f=(f_1,f_2)$, with $f_1(x)=[\![x\in\{4,6\}]\!]$ and $f_2(x)=[\![x=6]\!]$. Since $S_n^1=S_n^{f_1}$ and $S_n^2=S_n^{f_2}$, a central limit theorem for the 2-dimensional vector $(S_n^1,S_n^2)'$ is feasible provided that the quantities in (\ref{ide:CLT mu}) and (\ref{ide:CLT Sigma}) are computable and the $(2\times 2)$ matrix $\Sigma$ is positive definite. For this we differentiate the generating function $F(x,y_1,y_2)$ to obtain
\begin{eqnarray}
\label{ide:gf ESn1}
\sum_{n=1}^\infty E(S_n^1)\,x^n&=&\frac{\partial F}{\partial y_1}(x,1,1)=\frac{pqx^2}{(1-x)^2},\\ 
\sum_{n=1}^\infty E(S_n^2)\,x^n&=&\frac{\partial F}{\partial y_2}(x,1,1)=\frac{p^2q^2x^4}{(1-x)^2},\\ 
\sum_{n=1}^\infty E(S_n^1\cdot(S_n^1-1))\,x^n&=&\frac{\partial^2 F}{\partial y_1^2}(x,1,1)=\frac{2p^2q^2x^4}{(1-x)^3},\\  
\sum_{n=1}^\infty E(S_n^2\cdot(S_n^2-1))\,x^n&=&\frac{\partial^2 F}{\partial y_2^2}(x,1,1)=\frac{2p^3q^4(1-qx)x^7}{(1-x)^3},\\ 
\sum_{n=1}^\infty E(S_n^1\cdot S_n^2)\,x^n&=&\frac{\partial^2 F}{\partial y_1y_2}(x,1,1),\nonumber\\
\label{ide:gf ESn1Sn2}
&=&\frac{p^2q^2(1-q(q-p)x^2-px)x^4}{(1-x)^3}. 
\end{eqnarray}
The coefficients of each of these generating functions can be easily
extracted using (\ref{ide:coef (1-x)mm}). Furthermore, since for all random variables $X$ and $Y$ with finite second moment it applies that $\Var(X)=\E(X\cdot(X-1))-\E(X) \cdot(\E(X)-1)$ and that $\Cov(X,Y)=\E(X\cdot Y)-\E(X)\cdot\E(Y)$, one can deduce from (\ref{ide:coef (1-x)mm}) and (\ref{ide:gf ESn1})-(\ref{ide:gf ESn1Sn2}) the following asymptotic formulae as $n\to\infty$
\begin{eqnarray}
\label{ide:asymp ESn1}
\frac{\E(S_n^1)}{n}&=&pq+O\!\left(\frac{1}{n}\right),\\ 
\frac{\E(S_n^2)}{n}&=&p^2q^2+O\!\left(\frac{1}{n}\right),\\ 
\frac{\Var(S_n^1)}{n}&=&pq(1-3pq)+O\!\left(\frac{1}{n}\right),\\ 
\frac{\Var(S_n^2)}{n}&=&p^2q^2(1-11pq^2+13pq^3+6p^2q^2)+O\!\left(\frac{1}{n}\right),\\ 
\label{ide:asymp CovSn12}
\frac{\Cov(S_n^1,S_n^2)}{n}&=&\frac{p^2q^2}{2} (7p-pq-5+9q^2)
+O\!\left(\frac{1}{n}\right). 
\end{eqnarray}
These identities make explicit the terms in (\ref{ide:CLT mu}) and
(\ref{ide:CLT Sigma}). Furthermore, using symbolic algebra software
and replacing $q=1-p$ one can determine that
\[\det\Sigma=p^3(1-p)^3(1-5p+14p^2-25p^3+28p^4-16p^5+4p^6),\]
which is strictly positive for $0<p<1$. Therefore $((S_n^1,S_n^2)'-n\mu)/\sqrt{n}$ converges to a 2-dimensional centered normal random vector with variance-covariance matrix $\Sigma$, where $\mu$ and $\Sigma$ can be determined from
(\ref{ide:asymp ESn1})-(\ref{ide:asymp CovSn12}) as defined in
(\ref{ide:CLT mu}) and (\ref{ide:CLT Sigma}). For instance, if
$p=q=1/2$ then
\[\mu=16\cdot\left[\begin{array}{c}
4\\
1
\end{array}\right]\qquad;\qquad
\Sigma=256\cdot\left[\begin{array}{cc}
16 & 4\\
4 & 13
\end{array}\right].\]

\section{Application to correlated modular patterns}
\label{sec:app2}

Correlated, modular patterns are important in the analysis of RNA
motifs. Many functional molecules can be represented by sequence
motifs made up of modules separated by relatively unconstrained
spacer sequences \cite{bour99}. This modularity implies that there are
many more chances to match the pattern within a longer sequence than
would be possible for a simple pattern or moderate-size compound
pattern~\cite{sabe97,knig03,knig05}. This fact can greatly alter
estimates of the statistical significance of matching such a pattern.
The correlations between modules primarily take the form of base
pairs, which are essential for bringing the parts of the active site
into the structural juxtaposition required for function.

More generally, in many kinds of biological sequence analysis one is
interested in patterns that include correlations or gaps.  We use
numbers to denote correlations.  For example, in the case of the
binary alphabet $\{a,b\}$,
\[1a2a2b1=\{aaaaaba, aababba, baaaabb, bababbb\},\]
where either $a$ or $b$ could appear in the positions marked $1$ and $2$.

A gap of length exactly $k$ is denoted as $\#_k$ whereas a gap of
length at least $k$ is denoted $\#_k...$. The symbol $\#$ is used as a
shorthand for $\#_1$; in particular, $\#_k=\#\cdots\#$ $k$ times.  If
the symbols $\#$, $\#_k$ or $\#_{\ge k}$ appear more than once in the
same pattern each appearance is independendent.  For instance
\[1a\#_2b1=\{aaaaba,aaabba,aababa,aabbba,baaabb,baabbb,
  bababb,babbbb\}.\]
Finally, $ab\#...baa\#_2...bb$ is the set of all strings of the form
$abxbaaybb$ where $x,y\in\{a,b\}^*$ are such that $|x|\ge1$ and
$|y|\ge2$. This pattern consists of an infinite number of strings. In
this case we refer to $ab$, $baa$ and $bb$ as the modules of the pattern.

As in the previous section, for this example we consider the binary
alphabet $\A=\{a,b\}$ and let $X=(X_n)_{n\ge1}$ be a sequence of
i.i.d.  $\A$-valued random variables with initial distribution
$P(X_1=a)=p$ and $P(X_1=b)=q$, with $p\cdot q>0$ and $p+q=1$.  We
study the number of non-overlapping occurrences of the pattern
$aa\#...ba$ in $X_1...X_n$, and the sooner-time of the pattern
$1a\#...b1$ in $X$.

\begin{figure}[t]\centering
  \includegraphics[height=6cm]{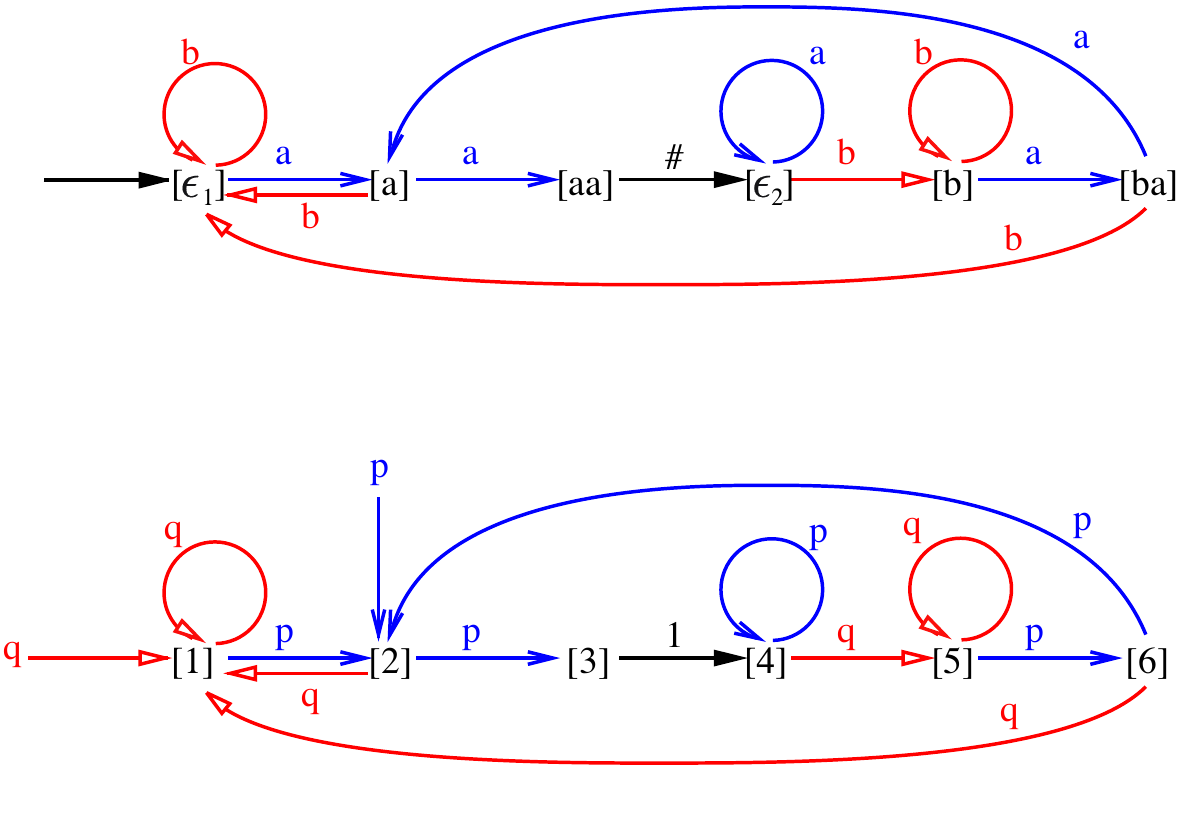}
 \caption[]{The automaton that recognizes the occurrence of the modular pattern
   $aa\#...ba$ in a text constructed using the binary alphabet $\{a,b\}$. Top,
   automaton $NC(aa\#...ba)$ which counts non-overlapping occurrences of
   $aa\#...ba$. The initial state is $\epsilon_1$ (left), and the
   terminal state is $ba$ (right). The symbol $\#$ corresponds to any
   of the two alphabet characters.  Bottom, representation of the first-order homogeneous Markov chain associated with a random text embedded in the automaton on top. The text is produced by a memoryless source where the character $a$ occurs with probability $p$ and the character $b$ occurs with probability $q$. The Markov chain starts at state $1$ with probability $q$ and at state $2$ with probability $p$.
   The probability that there are $m$ non-overlapping occurrences of $aa\#...ba$ in a random text of length $n$ is equivalent to the probability that the Markov chain visits state 6 a total $m$ times in the first $(n-1)$ steps. }
\label{fig:aastarba} 
\end{figure} 

\subsection{Frequency statistics of a modular pattern}
\label{sec:freq mod pat}

In addition to modular patterns that contain correlations through base
pairing, modular patterns that do not (as far as is currently known)
require base pairing are also important for processes involving RNA.
For example, transcriptional regulation requires combinatorial
regulation of binding sites for transcription factors that activate
and repress genes; splicing regulation requires specific combinations
of splicing enhancers and repression; and microRNA targeting appears
to be combinatorial \cite{sing06}. Many existing software packages for
detecting overrepresented words, such as the MobyDick package
\cite{buss00}, identify words that are surprisingly common given the
partition function by which they could be comprised of shorter words,
but fail to take into account correlations between word abundances
that could be caused by partial overlap of words of the same length.

To detect non-overlapping occurrences of the pattern $aa\#...ba$ in a
general text we first seek an automaton that detects the language
$\A^*aa\A\A^*ba$. This can be accomplished by concatenating the
Aho-Corasick automata $AC(\{aa\})$ and $AC(\{ba\})$: we concatenate the terminal state of $AC(\{aa\})$ with the
initial state of $AC(\{ba\})$ with two edges, one labeled with the
character $a$, and the other with the character $b$ (which we
represent visually as a single edge labeled with the character $\#$).
The resulting automaton is denoted as $AC(aa\#...ba)$. By definition,
the initial and terminal state of $AC(aa\#...ba)$ are the initial
state of $AC(\{ab\})$ and the terminal state of $AC(\{ba\})$,
respectively. The fact that $AC(aa\#...ba)$
recognizes the language $\A^*aa\A\A^*ba$ follows from
Theorem~\ref{thm:aho-corasick}.  Let $\epsilon_1$ and $\epsilon_2$
denote the initial state of $AC(\{aa\})$ and $AC(\{ba\})$
respectively.

To detect each non-overlapping occurrence of $aa\#...ba$ in the random
string $X_1...X_n$ we convey into the terminal state of
$AC(aa\#...ba)$ the transitions of its initial state as follows: first
remove all edges coming out from $ba$, add an edge labeled with the
character $a$ from $ba$ to $a$, and also add an edge labeled with the
character $b$ from $ba$ to $\epsilon_1$. We refer to this automaton as
$NC(aa\#...ba)$, or $NC$ in short. Here $N$ stands for non-overlapping
and $C$ for counting. See figure \ref{fig:aastarba} for a visual
representation of this automaton.

The fact that $NC(aa\#...ba)$ detects each
non-overlapping occurrence of the pattern $aa\#...ba$ in a general
text follows from the correctness of $AC(\{aa\})$ and $AC(\{ba\})$ and
the way these two automata were concatenated.

According to Theorem~\ref{thm:XG}, $X^{NC}$ is a first-order
homogeneous Markov chain with states $\epsilon_1$, $a$, $aa$,
$\epsilon_2$, $b$ and $ba$, which we label respectively as $1$, $2$,
$3$, $4$, $5$ and $6$. The initial distribution and probability
transition matrix of $X^{NC}$ are 
\[\mu:=[\begin{array}{cccccc}
q & p & 0 & 0 & 0 & 0
\end{array}]
\qquad;\qquad
P:=\left[\begin{array}{cccccc}
q & p & 0 & 0 & 0 & 0 \\
q & 0 & p & 0 & 0 & 0\\
0 & 0 & 0 & 1 & 0 & 0\\
0 & 0 & 0 & p & q & 0\\
0 & 0 & 0 & 0 & q & p\\
q & p & 0 & 0 & 0 & 0
\end{array}\right].\]
A visual representation of this Markov chain is displayed in figure
\ref{fig:aastarba}.

For $n\ge1$, we define
\[S_n:=\left(\begin{array}{c}
\hbox{number of non-overlapping occurrences}\\
\hbox{of $aa\#...ba$ as a substring of $X_1...X_n$}
\end{array}\right).\]
The distribution of $S_n$ corresponds to the number of visits that
$(X_i^{NC})_{i=1..n}$ makes to state $6$. The transfer matrix method
used in section \ref{sec:freq non-reduced} can now be used to
characterize the distribution of $S_n$. We therefore mark the edges
that are incident to state $6$ with a dummy variable $y$ that keeps
track of the number of times that this state is visited. This is equivalent to considering the matrix and vector
\[P_y:=\left[\begin{array}{cccccc}
q & p & 0 & 0 & 0 & 0 \\
q & 0 & p & 0 & 0 & 0\\
0 & 0 & 0 & 1 & 0 & 0\\
0 & 0 & 0 & p & q & 0\\
0 & 0 & 0 & 0 & q & py\\
q & p & 0 & 0 & 0 & 0
\end{array}\right]\qquad;\qquad
\delta:=\left[\begin{array}{c}
1\\
1\\
1\\
1\\
1\\
1
\end{array}\right].\]
As in (\ref{ide:prob S_n12=m1m2}), it follows
that
\[\Prob(S_n=m)=[y^m](\mu\cdot  P_y^{n-1}\cdot\delta)\qquad(n\ge1;m\ge0).\]
This result suffices to determine the exact distribution of
$S_n$ for small values of $n$. Furthermore, as in (\ref{ide:gen fct
  Sn12}), the generating function associated with $S_n$ is 
\begin{eqnarray*}
F(x,y)&:=&\sum_{n=1}^\infty\sum_{m=0}^\infty P(S_n=m)x^ny^m,\\ 
&=&x\cdot\mu\cdot(\Id_6-x\cdot P_y)^{-1}\cdot\delta,\\
&=&\frac{x(p^3qyx^4+p^2qx^3-qx+1)}{ 1-(p+2q)x+qx^2+p^2qx^3-
  p^2q^2x^4-p^3qyx^5}, 
\end{eqnarray*}
where the last identity was determined by using symbolic algebra
software to invert the matrix $(\Id_6-x\cdot P_y)$.

As in section~\ref{sec:freq non-reduced}, asymptotic formulae for
$\E(S_n)$ and $\Var(S_n)$ can be obtained using the partial fraction
decomposition of $\frac{\partial F}{\partial y}(x,1)$ and
$\frac{\partial^2 F}{\partial y^2}(x,1)$.  Since $X^{NC}$ is
irreducible and aperiodic,
\[\frac{S_n-\E(S_n)}{\sqrt{\Var(S_n)}}\]
can be shown to converge to a standard normal distribution.

\begin{figure}[t]\centering
  \includegraphics[height=10cm]{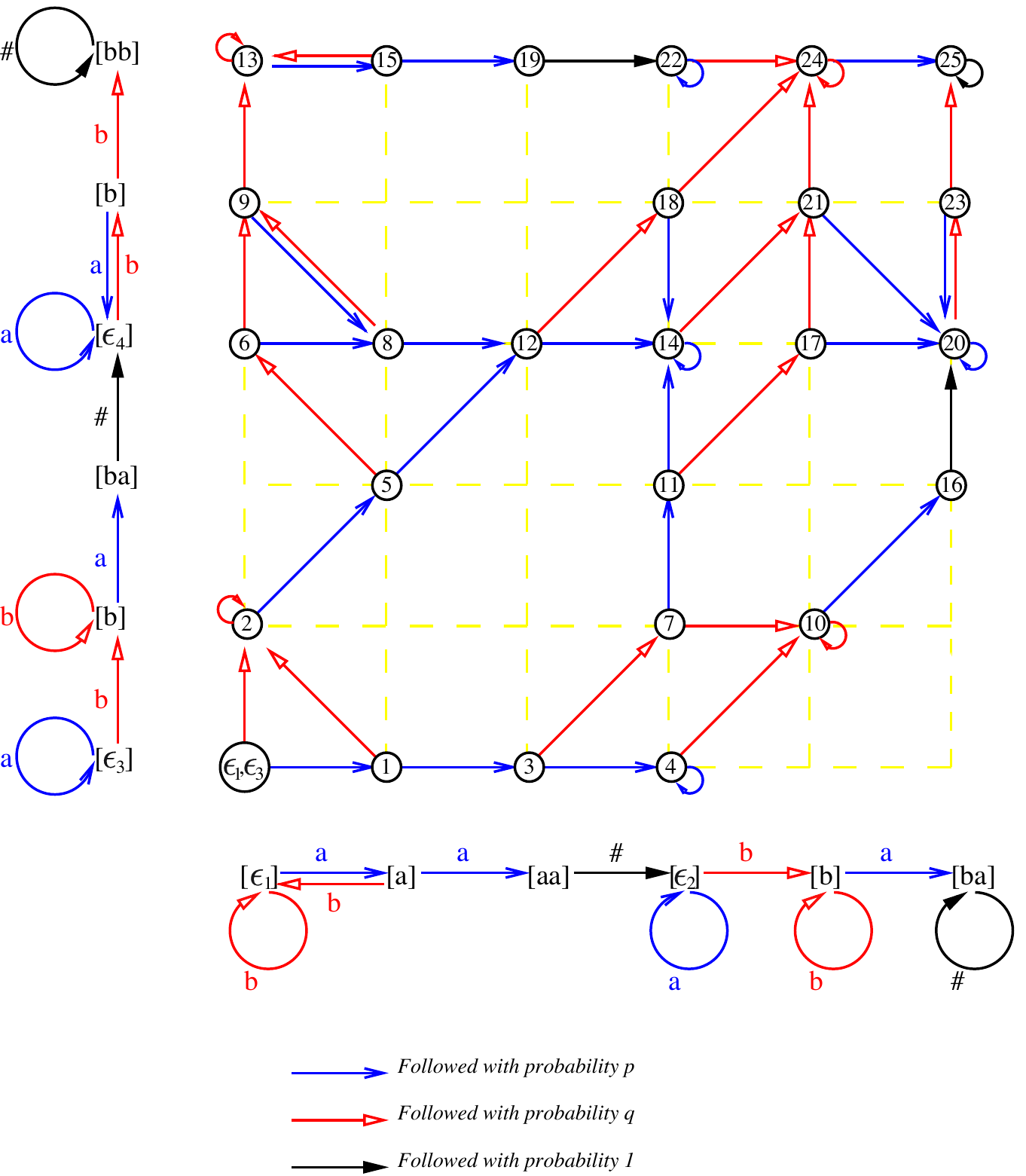}
 \caption[]{The automaton that recognizes the correlated modular pattern
$1a\#...b1$ in a text constructed using the binary alphabet $\{a,b\}$, where the symbol $1$ can be either $a$ or $b$ but must be the same character in both
   places. The grid in the middle is a visual representation of the
   synchronized automaton $ST(1a\#...b1)$, where horizontal axis is the automaton $ST(aa\#...ba)$ and the vertical axis is the automaton $ST(ba\#...bb)$. In the probabilistic automaton, blue edges are followed with
   probability $p$, red edges are followed with probability $q$, and
   black edges are followed with probability $1$. The terminal states
   13, 15, 19, 22, and 24 (top row) correspond to occurrence of the
   pattern $ba\#...bb$ but not $aa\#...ba$, while terminal states 16, 20
   and 23 (right column) correspond to occurrence of the pattern
   $aa\#...ba$ but not $ba\#...bb$. The terminal state 25 (top right
   corner) corresponds to recognition of both patterns.}
\label{fig:onebstara} 
\end{figure}

\subsection{Sooner-time of a correlated modular pattern}

Many functional RNAs must occur in a specific sequence context in
order to function. For example, riboswitches (RNA molecules that
regulate certain genes) must appear immediately upstream from the
start of the coding sequence \cite{wink02}. A related example is IRE,
the iron-responsive element in the ferritin mRNA, which binds a
protein cofactor to enhance the transcription of genes involved in
iron metabolism \cite{hent87}. In these cases, we are interested in
the distribution of first occurrences of a modular RNA pattern,
including correlations, relative to a specified start site.

To illustrate calculations of this type, define
\[T:=\hbox{ sooner-time of $1a\#...b1$ in $X$}.\]
In other words, $T$ is the smallest $n$ such that $X_1...X_n=...aa\#...ba$ or $X_1...X_n=...ba\#...bb$. To study the
distribution of $T$, we synchronize any automata that recognize the
languages $\A^*aa\A\A^*ba$ and $\A^*ba\A\A^*bb$. Consider the
automaton $AC(aa\#...ba)$ as defined in section \ref{sec:freq mod pat}.
 Similarly define $AC(ba\#...bb)$. Since
we are only interested in the number of occurrences of either of the
patterns $aa\#...ba$ or $ba\#...bb$, we turn the terminal states of
these automata into absorbing states. This is accomplished by
resetting all the edges coming out from terminal states to point to
themselves. We refer to the resulting automata as $ST(aa\#...ba)$ and
$ST(ba\#...bb)$ respectively, where $ST$ is short for sooner-time. A
visual representation of these automata can be 
found in figure \ref{fig:onebstara}.

Define $ST(1a\#...b1)$ to be the product of the automaton
$ST(aa\#...ba)$ with $ST(ba\#...bb)$. In principle, $ST(1a\#...b1)$
has 36 states, however, only 25 of these are accessible from the
initial state. For example, according to
Lemma~\ref{lem:aho-corasick}, there is no string
$x\in\A^*$ such that the path associated with $x$ in $ST(1a\#...b1)$
ends at state $(aa,ba)$.
A visual representation of $ST(1a\#...b1)$ reduced to only those
states that are accessible from the initial state is displayed in the
middle grid in figure \ref{fig:onebstara}, where for convenience we have relabeled
the accessible states as $1,\ldots,25$.

According to Theorem~\ref{thm:synchronization}, the pattern
$1a\#...b1$ occurs in a text provided that the path associated with
the text in $ST(1a\#...b1)$ ends at any of the states $13$, $15$,
$16$, $19$, $20$, $22$, $23$, $24$ and $25$. Furthermore, the
sooner-time of $1a\#...b1$ corresponds to the first time that any of
these states is visited.

To characterize the distribution of $T$, consider the first-order
homogeneous Markov chain $X^{ST}$. We denote the initial distribution
and probability transition matrix of $X^{ST}$ respectively as $\mu$
and $P$. Here $\mu$ is a row vector of dimension 25. The matrix $P$
has dimensions $25\times 25$ but is sparse (in each row there are only
two non-zero elements). By means of $X^{ST}$ the distribution of $T$
can be determined as shown in section \ref{sec:sooner-time non red}
for the sooner-time of a pair of non-reduced patterns. Define

\[
\nu:=\left[\begin{array}{c}
p\\
q\\
0\\
0\\
0\\
0\\
0\\
0\\
0\\
0\\
0\\
0\\
0\\
0\\
0\\
0
\end{array}\right]
\,;\,
Q:=\left[\begin{array}{ccccccccccccccccc}
0 & q & p & 0 & 0 & 0 & 0 & 0 & 0 & 0 & 0 & 0 & 0 & 0 & 0 & 0 \\
0 & q & 0 & 0 & p & 0 & 0 & 0 & 0 & 0 & 0 & 0 & 0 & 0 & 0 & 0 \\
0 & 0 & 0 & p & 0 & 0 & q & 0 & 0 & 0 & 0 & 0 & 0 & 0 & 0 & 0 \\
0 & 0 & 0 & p & 0 & 0 & 0 & 0 & 0 & q & 0 & 0 & 0 & 0 & 0 & 0 \\
0 & 0 & 0 & 0 & 0 & q & 0 & 0 & 0 & 0 & 0 & p & 0 & 0 & 0 & 0 \\
0 & 0 & 0 & 0 & 0 & 0 & 0 & p & q & 0 & 0 & 0 & 0 & 0 & 0 & 0 \\
0 & 0 & 0 & 0 & 0 & 0 & 0 & 0 & 0 & q & p & 0 & 0 & 0 & 0 & 0 \\
0 & 0 & 0 & 0 & 0 & 0 & 0 & 0 & q & 0 & 0 & p & 0 & 0 & 0 & 0 \\
0 & 0 & 0 & 0 & 0 & 0 & 0 & p & 0 & 0 & 0 & 0 & 0 & 0 & 0 & 0 \\
0 & 0 & 0 & 0 & 0 & 0 & 0 & 0 & 0 & q & 0 & 0 & 0 & 0 & 0 & 0 \\
0 & 0 & 0 & 0 & 0 & 0 & 0 & 0 & 0 & 0 & 0 & 0 & p & q & 0 & 0 \\
0 & 0 & 0 & 0 & 0 & 0 & 0 & 0 & 0 & 0 & 0 & 0 & p & 0 & q & 0 \\
0 & 0 & 0 & 0 & 0 & 0 & 0 & 0 & 0 & 0 & 0 & 0 & p & 0 & 0 & q \\
0 & 0 & 0 & 0 & 0 & 0 & 0 & 0 & 0 & 0 & 0 & 0 & 0 & 0 & 0 & q \\
0 & 0 & 0 & 0 & 0 & 0 & 0 & 0 & 0 & 0 & 0 & 0 & p & 0 & 0 & 0 \\
0 & 0 & 0 & 0 & 0 & 0 & 0 & 0 & 0 & 0 & 0 & 0 & 0 & 0 & 0 & 0 \\
\end{array}\right]
\,;\,
u:=\left[\begin{array}{c}
0\\
0\\
0\\
0\\
0\\
0\\
0\\
0\\
q\\
p\\
0\\
0\\
0\\
p\\
q\\
1
\end{array}\right].\]
The vector $\nu$ corresponds to the vector $\mu$ with columns $13$,
$15$, $16$, $19$, $20$, $22$, $23$, $24$ and $25$ removed. The matrix
$Q$ corresponds to the matrix $P$ but with rows and columns $13$,
$15$, $16$, $19$, $20$, $22$, $23$, $24$ and $25$ removed.  Finally,
the vector $u$ corresponds to the sum of columns $13$, $15$, $16$,
$19$, $20$, $22$, $23$, $24$ and $25$ in $P$ however with the rows of
these same number removed. Since states $13$, $15$, $16$, $19$, $20$,
$22$, $23$, $24$ and $25$ correspond to the detection of the pattern
$1a\#...b1$, it follows that
\[\Prob(T=n)=\nu\cdot Q^{n-2}\cdot u\qquad(n\ge2).\]
The generating function associated with $T$ is 
\begin{eqnarray*}
F(x)&:=&\sum_{n=2}^\infty \Prob(T=n)x^n,\\
&=&x^2\cdot\nu\cdot(\Id-x\cdot Q)^{-1}\cdot u,\\
&=&\frac{pqx^5(2p^3q^3x^5+p^3q^3x^4-3p^3q^3x^3
  +p^2q^2x^2+pq(1-pq)x+p^2+q^2)}{(1-px)(1-qx)(1-pqx^2)}, 
\end{eqnarray*}
where for the last identity we have used 
symbolic algebra software.

This result for the generating function provides useful information
about the distribution of $T$.  For instance, if $p\ne q$ then
$x=\min\{1/p,1/q\}$ is a simple zero and the closest zero to the
origin of the denominator $F(x)$. On the other hand, since
$(p^2q^2x^2-3p^3q^3x^3)\ge0$ for all $x\in[0,1]$, the numerator of
$F(x)$ does not vanish at $x=\min\{1/p,1/q\}$. Hence $\Prob[T=n]\sim
c_1(p,q)\cdot (\min\{p,q\})^{-n}$ as $n\to\infty$, where the constant
$c_1(p,q)>0$ is a computable constant that can be determined from the
partial fraction decomposition of $F(x)$.

For the case $p=q$, we find that
\[F(x)=\frac{x^5(2x^4-3x^3+3x^2-2x+16)}{16(2-x)^3}.\]
In this case, $x=2$ is a zero (of order 3) of the denominator of
$F(x)$ but not of its numerator. Using the partial fraction
decomposition of $F(x)$ and (\ref{ide:coef (1-x)mm}), it follows that
$\Prob[T=n]\sim c_2\cdot n^2/2^n$ as $n\to\infty$, where $c_2$ is a computable constant from the partial fraction decomposition of $F(x)$. 

Refinements of this argument can be used to find explicit formulae for
the probabilities and generating functions associated with the events
$X_1...X_T=...aa\#...ba$ and $X_1...X_T=...ba\#...bb$.

\section{Conclusions}

In this paper we have reviewed the use of
deterministic finite automata for probabilistic pattern matching. This
view of the pattern matching problem allows many different problems to
be addressed in a general framework, and unifies different
ideas addressed in the computer science, mathematics, and
bioinformatics literature. We have summarized the key results
to present a self-contained mathematical summary of
previous work, including definitions, theorems, proofs, and examples.

The key results of deterministic automata are how to construct state
machines from possibly simpler state machines
to find matches of regular patterns in a given text. The Aho-Corasick automaton (based on the maximum prefix-suffix rule) is  a classic example of an automaton that
recognizes a set of keywords in a text.
For matching compound patterns (i.e., containing multiple keywords), the
synchronization of multiple automata is an important tool. 

For assessing the statistical significance of motif searches in
biological sequence data, the pattern matching problem must be
extended to determine the probability that a pattern occurs in a
random string (specified by a given model). Mathematically, this can
be done by the Markov chain embedding of a random text into an automaton. This means
considering a random walk on the automaton where the transition
probabilities of the walk are determined by the model which
generates the random string. This maps the probabilistic pattern
matching problem onto a Markov chain, allowing techniques
from combinatorics and the theory of Markov chains to be applied to the problem. In
particular, the probability that a given set of patterns occurs in the
random text corresponds to the probability that a specific Markov
chain visits a certain set of terminal states.

To illustrate the application of these ideas to biological sequence
analysis, we presented two examples. In all the examples, we used a
simplified binary alphabet and patterns that admit a simple description to illustrate the
key ideas.

The first application was the search for a compound pattern consisting of
two keywords. We
demonstrated how to determine the transition matrix of the Markov
chain which determines the probability that one of the keywords occurs in the random string. We then derived the sooner-time probability
distribution, the probability that any of the keywords first
occurs after $n$ characters in the random string, and used generating
function methods to derive the asymptotic distributions for large $n$.
We used similar mathematical methods to derive the probability that either of the two keywords occurs a given number of times in a random string of $n$ characters.

The second application was the search for a correlated modular
pattern, in which two sub-patterns (modules) must appear in a certain
order but can be separated by an arbitrary number of characters.
Correlations mean that certain characters within the pattern can take
different values, but the values must be correlated (for example,
through base pairing).  We illustrated the calculation of the
frequency statistic of a modular pattern, including the asymptotic
probability distribution for large $n$. We also derived
formulae for the probability that a correlated modular pattern first
appears after $n$ characters in the string.

We have focused in this review on random strings produced by
memoryless sources. However, we have provided references for Markovian sources and hidden Markov models (which can also be handled in this
framework).

These methods are applicable to determining the significance of motif
searches in genome sequences. In particular, modular and correlated
patterns frequently occur in the sequences of functional RNA
molecules.  As the number of functional RNAs increases, the ability to
infer the statistical significance of
matches to RNA sequence patterns is increasingly important.  A unified mathematical framework (based on the concepts of
automata, Markov chain embedding, and synchronization) can be used to
analyze a range of biologically important pattern-matching problems,
including the regulation of splicing and transcription and the
probability of occurrence of catalytic RNA motifs in genomes or
random-sequence RNA pools.  All these apparently different problems can
be addressed in the framework of probabilistic pattern matching.

\bibliographystyle{alpha}
\bibliography{biblio}

\end{document}